\documentclass[11pt]{article}
\usepackage{fullpage}

\usepackage{epsf}
\usepackage{fancyhdr}
\usepackage{graphics}
\usepackage{graphicx}
\usepackage{float}
\usepackage{enumerate}

\usepackage{color}

\usepackage{amsthm}
\usepackage{amsmath}
\usepackage{amssymb}

\usepackage[numbers]{natbib}
\usepackage{tikz}
\usepackage{hyperref}
\usepackage{mathtools}
\usepackage{mathrsfs}
\mathtoolsset{showonlyrefs}

\theoremstyle{plain}

\newtheorem{theo}{Theorem}[section]

\newtheorem{lem}{Lemma}[section]
\newtheorem{prop}{Proposition}[section]
\newtheorem{cor}{Corollary}[section]

\theoremstyle{definition} 

\newtheorem{nota}{Notation}[section]
\newtheorem{de}{Definition}[section]
\newtheorem{exa}{Example}[section]
\newtheorem{as}{Assumption}[section]
\newtheorem{alg}{Algorithm}[section]

\newcommand{\btheo}{\begin{theo}}
\newcommand{\bde}{\begin{de}}
\newcommand{\ble}{\begin{lem}}
\newcommand{\bpr}{\begin{prop}}
\newcommand{\bno}{\begin{nota}}
\newcommand{\bex}{\begin{exa}}
\newcommand{\bcor}{\begin{cor}}
\newcommand{\spro}{\begin{proof}}
\newcommand{\bas}{\begin{as}}
\newcommand{\balg}{\begin{alg}}

\newcommand{\etheo}{\end{theo}}
\newcommand{\ede}{\end{de}}
\newcommand{\ele}{\end{lem}}
\newcommand{\epr}{\end{prop}}
\newcommand{\eno}{\end{nota}}
\newcommand{\eex}{\end{exa}}
\newcommand{\ecor}{\end{cor}}
\newcommand{\fpro}{\end{proof}}
\newcommand{\eas}{\end{as}}
\newcommand{\ealg}{\end{alg}}

\theoremstyle{plain}

\newtheorem{theos}{Theorem}
\newtheorem{props}{Proposition}
\newtheorem{lems}{Lemma}
\newtheorem{cors}{Corollary}

\theoremstyle{definition}
\newtheorem{exas}{Example}
\newtheorem{algs}{Algorithm}
\newtheorem{asss}{Assumption}
\newtheorem{defns}{Definition}

\newcommand{\btheos}{\begin{theos}}
\newcommand{\etheos}{\end{theos}}
\newcommand{\brems}{\begin{remark}}
\newcommand{\erems}{\end{remark}}
\newcommand{\bprops}{\begin{props}}
\newcommand{\eprops}{\end{props}}
\newcommand{\bdes}{\begin{defns}}
\newcommand{\edes}{\end{defns}}
\newcommand{\blems}{\begin{lems}}
\newcommand{\elems}{\end{lems}}
\newcommand{\bcors}{\begin{cors}}
\newcommand{\ecors}{\end{cors}}
\newcommand{\bexs}{\begin{exas}}
\newcommand{\eexs}{\end{exas}}
\newcommand{\balgs}{\begin{algs}}
\newcommand{\ealgs}{\end{algs}}
\newcommand{\bass}{\begin{asss}}
\newcommand{\eass}{\end{asss}}

\newcommand{\bigrad}{\ensuremath{R}}

\newcommand{\myspan}{\ensuremath{\operatorname{span}}}

\newcommand{\bcar}{\begin{itemize}}
\newcommand{\ecar}{\end{itemize}}

\newcommand{\Constraint}{\ensuremath{\mathcal{S}}}

\newcommand{\mytrace}[1]{\ensuremath{\operatorname{tr} \big(#1 \big)}}

\newcommand{\Prob}{\ensuremath{\mprob}}

\newcommand{\Vset}{\ensuremath{\mathcal{V}}}

\newcommand{\Aevent}{\ensuremath{\mathcal{A}}}

\newcommand{\smoothpar}{\ensuremath{\alpha}}

\newcommand{\CritFun}{\ensuremath{F}}

\newcommand{\Mpack}{\ensuremath{M}}

\newcommand{\MiniMax}{\ensuremath{\mathfrak{M}}}

\newcommand{\zhat}{\ensuremath{\widehat{z}}}

\newcommand{\genparam}{\ensuremath{\theta}}
\newcommand{\genparamstar}{\ensuremath{\genparam^*}}
\newcommand{\genparamhat}{\ensuremath{\widehat{\genparam}}}

\newcommand{\thetastar}{\ensuremath{{\theta^*}}}

\newcommand{\ytil}{\ensuremath{\widetilde{y}}}

\newcommand{\NORMAL}{\ensuremath{\mathcal{N}}}

\newcommand{\diag}{\ensuremath{\operatorname{diag}}}

\newcommand{\kind}{\ensuremath{k}}

\newcommand{\diam}{\ensuremath{\operatorname{diam}}}

\newcommand{\thetahat}{\ensuremath{\widehat{\theta}}}

\newcommand{\usedim}{\ensuremath{d}}

\newcommand{\Ellipse}{\ensuremath{\mathcal{E}}}

\definecolor{MyGray}{rgb}{0.9,0.9,0.9}

\makeatletter\newenvironment{graybox}{
\begin{lrbox}{\@tempboxa}
\begin{minipage}{1\columnwidth}}{\end{minipage}
\end{lrbox}%
\colorbox{MyGray}{\usebox{\@tempboxa}} }
\makeatother



\newcommand{\Zback}[1]{\ensuremath{Z^{\backslash i}}}

\newcommand{\argmin}{\ensuremath{\operatorname{argmin}}}

\newcommand{\xsamstack}[1]{\ensuremath{x_1^\numobs}}
\newcommand{\Xsamstack}[1]{\ensuremath{X_1^\numobs}}

\newcommand{\gtil}{\ensuremath{\widetilde{g}}}

\newcommand{\widgraph}[2]{\includegraphics[keepaspectratio,width=#1]{#2}}

\newcommand{\noisestd}{\ensuremath{\sigma}}

\newcommand{\delcrit}{\ensuremath{\delta_\numobs}}

\newcommand{\Term}{\ensuremath{T}}

\newcommand{\Tset}{\ensuremath{T}}

\newcommand{\fancysoln}[1]{
\ifthenelse{\equal{\doctype}{WITHSOLS}}
{
\begin{soln}
#1
\end{soln}
}
{
}
}

\newcommand{\fstar}{\ensuremath{f^*}}

\long\def\comment#1{}

\makeatletter
\def\@cite#1#2{[\if@tempswa #2 \fi #1]}
\makeatother

\newcommand{\Hil}{\ensuremath{\mathcal{H}}}

\newcommand{\fhat}{\ensuremath{\widehat{f}}}

\newcommand{\Xspace}{\ensuremath{\mathcal{X}}}

\newcommand{\defn}{\ensuremath{: \, = }}

\newcommand{\var}{\ensuremath{\operatorname{var}}}

\newcommand{\real}{\ensuremath{\mathbb{R}}}

\newcommand{\pdim}{\ensuremath{p}}
\newcommand{\numobs}{\ensuremath{n}}

\newcommand{\Sset}{\ensuremath{S}}

\newcommand{\mprob}{\ensuremath{\mathbb{P}}}

\newcommand{\Ball}{\ensuremath{\mathbb{B}}}

\newcommand{\sgauss}{\ensuremath{w}}

\newcommand{\Amat}{\ensuremath{\mymathbf{A}}}

\newcommand{\betastar}{\ensuremath{\beta^*}}

\newcommand{\betahat}{{\ensuremath{{\widehat{\beta}}}}}

\newlength{\widebarargwidth}
\newlength{\widebarargheight}
\newlength{\widebarargdepth}

\newcommand{\vecnorm}[2]{\| #1\|_{#2}}

\newcommand{\inprod}[2]{\ensuremath{\langle #1 , \, #2 \rangle}}

\newcommand{\kullinf}[2]{\ensuremath{D(#1\, \| \, #2)}}

\newcommand{\Exs}{\ensuremath{\mathbb{E}}}

\newcommand{\eigplain}{\ensuremath{\mu}}
\newcommand{\eig}[1]{\ensuremath{\eigplain_{#1}}}

\newcommand{\deltachat}{\ensuremath{\delta_0}}

\newcommand{\Deltatilde}{\ensuremath{\widetilde{\Delta}}}
\newcommand{\deltatilde}{\ensuremath{\widetilde{\delta}}}
\renewcommand{\delcrit}{\ensuremath{\delta_*}}
\newcommand{\delcritleft}{\ensuremath{\delta'}}
\newcommand{\delcritright}{\ensuremath{\delta''}}

\newcommand{\CritFunPlain}{\ensuremath{g}}
\renewcommand{\CritFun}{\ensuremath{\CritFunPlain}}
\newcommand{\CritFunUp}{\ensuremath{\CritFun^u}}
\newcommand{\CritFunLw}{\ensuremath{\CritFun^\ell}}
\newcommand{\CritFunTilde}{\ensuremath{\widetilde{\CritFunPlain}}}

\newcommand{\enorm}[1]{\ensuremath{\vecnorm{#1}{\Ellipse}}}

\newcommand{\Proj}{\ensuremath{\Pi}}
\newcommand{\ProjClass}[1]{\ensuremath{\mathcal{P}_{#1}}}
\newcommand{\Indic}[1]{\ensuremath{\mathbf{1}}\{#1\}}
\DeclarePairedDelimiter{\floor}{\lfloor}{\rfloor}

\newcommand{\IdMat}{\ensuremath{\mathbf{I}}}
\renewcommand{\Amat}{\ensuremath{A}}

\newcommand{\thetadag}{\ensuremath{\theta^\dagger}}
\newcommand{\at}{\genparamstar}

\newcommand{\ltwo}[1]{\ensuremath{\|#1\|_2}}
\newcommand{\GWidth}{\ensuremath{\mathscr{G}}}
\newcommand{\kwidth}{\ensuremath{\mathscr{W}}}

\renewcommand{\kind}{\ensuremath{k_*}}
\newcommand{\Sph}{\ensuremath{\mathbb{S}}}
\renewcommand{\Ball}{\ensuremath{\mathbb{B}}}

\newcommand{\upk}{\ensuremath{\kind}}
\newcommand{\upktilde}{\ensuremath{\widetilde{\kind}}}
\newcommand{\lwk}{\ensuremath{\kind}}

\renewcommand{\Aevent}{\ensuremath{A}}

\newcommand{\Ncover}{\ensuremath{N}}
\newcommand{\trunc}[1]{\Pi_{#1}}
\newcommand{\LPset}[1]{\ensuremath{\mathcal{P}_{#1}}}
\newcommand{\Mat}{H}
\newcommand{\NewMat}{M}

\DeclarePairedDelimiter{\parens}{(}{)}
\DeclarePairedDelimiter{\braces}{\{}{\}}

\newcommand{\ccon}{\ensuremath{c}}

\newcommand{\newepscritu}{t^*_u}
\newcommand{\newepscritl}{t^*_{\ell}}

\newcommand{\ellip}{\mathcal{E}}

\renewcommand{\Tset}{\ensuremath{T}}
\renewcommand{\Sset}{\ensuremath{S}}
\newcommand{\Gammaclass}{\ensuremath{\Gamma}}
\newcommand{\zs}{\ensuremath{z^{\Sset}}}

\newcommand{\Rfun}{\ensuremath{\Phi}}
\newcommand{\InvRfun}{\Rfun^{-1}}

\newcommand{\Tclass}{\mathscr{T}}
\newcommand{\gring}{\ensuremath{\widehat{g}}}

\renewcommand{\Aevent}{\ensuremath{\mathcal{A}}}

\newcommand{\tinyconst}{\ensuremath{\eta}}

\newcommand{\constlw}{\ensuremath{c_\ell}}
\newcommand{\constup}{\ensuremath{c_u}}

\newcommand{\indic}{\ensuremath{\mathbf{1}}}

\newcommand{\KerPlain}{\ensuremath{\mathcal{K}}}
\newcommand{\PlainKer}{\KerPlain}

\makeatletter
\long\def\@makecaption#1#2{
        \vskip 0.8ex
        \setbox\@tempboxa\hbox{\small {\bf #1:} #2}
        \parindent 1.5em  
        \dimen0=\hsize
        \advance\dimen0 by -3em
        \ifdim \wd\@tempboxa >\dimen0
                \hbox to \hsize{
                        \parindent 0em
                        \hfil
                        \parbox{\dimen0}{\def\baselinestretch{0.96}\small
                                {\bf #1.} #2
                                }
                        \hfil}
        \else \hbox to \hsize{\hfil \box\@tempboxa \hfil}
        \fi
        }
\makeatother

\begin{document}


\begin{center}

{\bf{\LARGE{From Gauss to Kolmogorov:\\Localized Measures of
      Complexity for Ellipses}}}

\vspace*{.2in}

\large{
\begin{tabular}{ccc}
Yuting Wei$^\dagger$ & Billy Fang$^\dagger$ & Martin
J. Wainwright$^{\dagger,\star}$
\end{tabular}
}

\vspace*{.2in}

\begin{tabular}{c}
Department of Statistics$^\dagger$, and \\ Department of Electrical
Engineering and Computer Sciences$^\star$ \\ UC Berkeley, Berkeley, CA
94720
\end{tabular}

\vspace*{.2in}

\today

\vspace*{.2in}

\begin{abstract}
The Gaussian width is a fundamental quantity in probability,
statistics and geometry, known to underlie the intrinsic difficulty of
estimation and hypothesis testing. In this work, we show how the Gaussian width,
when localized to any given point of an ellipse, can be controlled by
the Kolmogorov width of a set similarly localized.  This connection
leads to an explicit characterization of the estimation error of
least-squares regression as a function of the true regression vector
within the ellipse.  The rate of error decay varies substantially as a
function of location: as a concrete example, in Sobolev ellipses of
smoothness $\alpha$, we exhibit rates that vary from
$(\sigma^2)^{\frac{2 \alpha}{2 \alpha + 1}}$, corresponding to the
classical global rate, to the faster rate $(\sigma^2)^{\frac{4
    \alpha}{4 \alpha + 1}}$.  We also show how the local Kolmogorov
width can be related to local metric entropy.
\end{abstract}
\end{center}



\section{Introduction}
\label{SecIntroduction}

The Gaussian width is an important measure of the complexity of a set,
and it plays an important role in geometry, statistics and probability
theory.  Most relevant to this paper is its central role in empirical
process theory, where the Gaussian width and its Bernoulli analogue
(known as the Rademacher width) can be used to upper bound the error
for various types of non-parametric
estimators~\cite{vandeGeer,vanderVaart96,Bar05,
  koltchinskii2001rademacher,chatterjee2014new}.  More recently, these
same complexity measures have also been shown to play an important
role in high-dimensional testing problems~\cite{wei2017geometry}.

For a general set, it is non-trivial to provide analytical expressions
for its Gaussian or Rademacher widths.  There are a variety of
techniques for obtaining bounds, including upper bounds via the
classical entropy integral of Dudley, as well as lower bounds due to
Sudakov-Fernique (see the book~\cite{LedTal91} for details on these
and other results).  More recently, \citet{Tal00} has introduced a
generic chaining technique that leads to sharp lower and upper bounds.
However, for a general set, it is impossible to evaluate the
expressions obtained from the generic chaining, and so for
applications in statistics, it is of considerable interest to develop
techniques that yield tractable characterizations of various forms of
widths.

In this paper, we study a class of Gaussian widths that arise in the
context of estimation over (possibly infinite-dimensional) ellipses.
As we describe below, many non-parametric problems, among them
are regression and density estimation over classes of smooth functions,
can be reduced to such ellipse estimation problems.  Obtaining sharp
rates for such estimation problems requires studying a
\emph{localized} notion of Gaussian width, in which the ellipse is
intersected with a Euclidean ball around the element $\thetastar$
being estimated.  The main technical contribution of this paper is to
show how this localized Gaussian width can be bounded, from both above
and below, using a localized form of the Kolmogorov
width~\cite{pinkus2012n}.  As we show with a number of corollaries,
this Kolmogorov width can be calculated in many interesting examples.

Our work makes a connection to the evolving line of work on
instance-specific rates in estimation and testing.  Within the
decision-theoretic framework, the classical approach is to study the
(global) minimax risk over a certain problem class.  In this
framework, methods are compared via their worst-case behavior as
measured by performance over the entire problem class.  For the
ellipse problems considered here, global minimax risks in  various
norms are well-understood; for instance, see the classic
papers~\cite{pinsker1980optimal,ibragimov1978,ibragimov2013statistical}.
When the risk function is near to constant over the set, then the
global minimax risk is reflective of the typical behavior.  If not,
then one is motivated to seek more refined ways of characterizing the
hardness of different problems, and the performance of different
estimators.

One way of doing so is by studying the notion of an adaptive
estimator, meaning one whose performance automatically adapts to some
(unknown) property of the underlying function being estimated.  For
instance, estimators using wavelet bases are known to be adaptive to
unknown degree of smoothness~\cite{Donoho95c,donoho1994ideal}.
Similarly, in the context of shape-constrained problems, there is a
line of work showing that for functions with simpler structure, it is
possible to achieve faster rates than the global minimax ones
(e.g. \cite{meyer2000degrees,zhang2002risk,chatterjee2015risk}).  A
related line of work, including some of our own, has studied
adaptivity in the context of hypothesis testing
(e.g.,~\cite{valiant2014automatic,
  balakrishnan2017hypothesis,wei2017testing}).  The adaptive
estimation rates established in this work also share this spirit of
being instance-specific.


\subsection{Some motivating examples}

A primary motivation for our work is to understand the behavior of
least-squares estimators over ellipses.  Accordingly, let us give a
precise definition of the ellipse estimation problem, along with some
motivating examples.

Given a fixed integer $\usedim$ and a sequence of non-negative scalars
$\eig{1} \geq \ldots \geq \eig{\usedim} \geq 0$, we can define an
elliptical norm on $\real^\usedim$ via $\enorm{\genparam}^2 \defn
\sum_{j=1}^\usedim \frac{\genparam_j^2}{\eig{j}}$.  Here for any
coefficient $\mu_k = 0$, we interpret the constraint as enforcing that
$\genparam_k= 0$.  For any radius $\bigrad > 0$, this semi-norm
defines an ellipse of the form
\begin{align}
\label{EqnEllipse}
\Ellipse(\bigrad) \defn \Big\{ \theta \in \real^\usedim ~\mid~
\enorm{\genparam} \leq \bigrad \Big \}.
\end{align}
We frequently focus on the case $\bigrad = 1$, in which case we adopt
the shorthand notation $\Ellipse$ for the set $\Ellipse(1)$.  Whereas
equation~\eqref{EqnEllipse} defines a finite-dimensional ellipse, it
should be noted that our theory also applies to infinite-dimensional
ellipses for sequences $\{\mu_j\}_{j=1}^\infty$ that are summable.
Such results can be recovered by studying a truncated version of the
ellipse with finite dimension $\usedim$, and then taking suitable
limits.  In order to simplify the exposition, we develop our results
with finite $\usedim$, noting how they extend to infinite dimensions
after stating our results.

Suppose that for some unknown vector $\thetastar \in \Ellipse$, we
make noisy observations of the form
\begin{align}
\label{EqnGaussianModel}
y = \genparamstar + \noisestd \sgauss, \qquad \text{where } \sgauss
\sim \NORMAL(0, \IdMat_\usedim).
\end{align}
We assume that the ellipse $\Ellipse$ and noise standard deviation
$\noisestd$ is known.  The goal of ellipse estimation is to specify a
mapping $y \mapsto \thetahat(y)$ such that the associated Euclidean risk
$\Exs_y \ltwo{\thetahat(y) - \thetastar}^2$ is as small as possible.

\noindent Let us consider some concrete problems that can be reduced
to instances of ellipse estimation.

\begin{exas}[Linear prediction with correlated designs]
Suppose that we make observations from the standard linear model
\begin{align*}
\ytil & = X \betastar + \nu w,
\end{align*}
where $\ytil \in \real^\numobs$ is the response vector, $X \in
\real^{\numobs \times \pdim}$ is a (fixed, non-random) design matrix,
and $w \sim N(0, \IdMat_{\numobs})$ is noise.  Suppose moreover that
we know a priori that $\|\betastar\|_2 \leq R$ for some radius $R >
0$.  Alternatively, we can think of a condition of this form arising
implicitly when using estimators such as ridge regression.

\begin{figure}[h!]
  \begin{center}
    \widgraph{.6\textwidth}{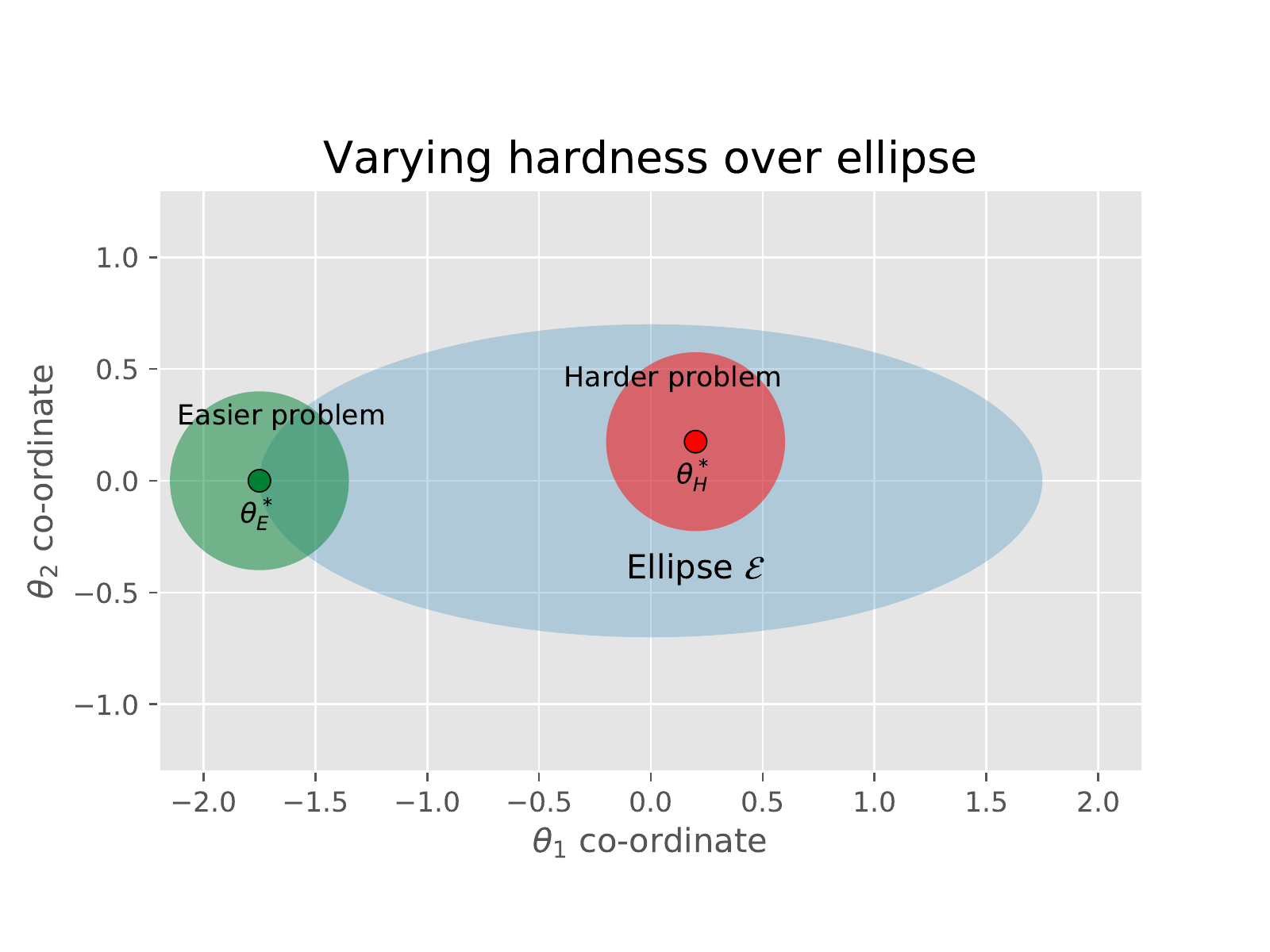}
    \caption{Illustration of the ellipse estimation problem.  The goal
      to estimate an unknown vector $\thetastar$ belonging to an
      ellipse based on noisy observations.  The local geometry of the
      ellipse controls the difficulty of the problem: due to its
      proximity to the narrow end of the ellipse, the vector
      $\theta^*_E$ is relatively easy to estimate.  By contrast, the
      vector $\theta^*_H$ should be harder, since it lies closest to
      the center of the ellipse.  The theory given in this paper
      confirms this intuition; see \autoref{SecEstimation} for
      details.}
    \label{FigHardEllipse}
  \end{center}
\end{figure}

Given an estimate $\betahat$, its prediction accuracy can be assessed
via the mean-squared error \mbox{$\Exs \big[ \frac{1}{\numobs} \|X
    \betahat - X \betastar\|_2^2 \big]$,} where the expectation is
taken over the observation noise.  Equivalently, letting $\thetahat =
X \betahat/\sqrt{\numobs}$ and $\thetastar = X
\betastar/\sqrt{\numobs}$, our problem is to minimize the mean-squared
error $\Exs \|\thetahat - \thetastar\|_2^2$.  After
this transformation, we arrive at the observation model $y =
\thetastar + \frac{\nu}{\sqrt{\numobs}} w$, which is a version of our
original model~\eqref{EqnGaussianModel} with $\usedim = \numobs$ and
$\sigma = \frac{\nu}{\sqrt{\numobs}}$.  Moreover, the constraint on
the $\ell_2$-norm of $\betastar$ translates into an ellipse constraint
on $\thetastar$.  In particular, the ellipse is determined by the
non-zero eigenvalues of the matrix $\frac{1}{\numobs} X X^\top \in \real^{\numobs \times
  \numobs}$.

As shown in \autoref{FigHardEllipse}, it is natural to conjecture
that the location of $\thetastar$ within this ellipse affects the
difficulty of estimation.  Note that $\Exs \|y - \thetastar\|_2^2 =
\nu^2 / \numobs$, so that on average, the observed vector $y$ lies at squared
Euclidean distance $\nu^2 / \numobs$ from the true vector.  In certain favorable
cases, such as a vector $\theta^*_E$ that lies at or close to the
boundary of an elongated side of the ellipse, the side-knowledge that
$\thetastar \in \Ellipse$ is helpful.  In other cases, such as a
vector $\theta^*_H$ that lies closer to the center of the ellipse, the
elliptical constraint is less helpful.  The theory to be developed in
this paper makes this intuition precise.  In particular,
\autoref{SecEstimation} is devoted to a number of consequences of our
main results for the problem of estimation in ellipses.

\end{exas}

\begin{exas}[Non-parametric regression using reproducing kernels]
  We now turn to a class of non-parametric problems that involve a
  form of ellipse estimation.  Suppose that our goal is to predict a
  response $z \in \real$ based on observing a collection of predictors
  $x \in \Xspace$.  Assuming that pairs $(X,Z)$ are drawn jointly from
  some unknown distribution $\mprob$, the optimal prediction in terms
  of mean-squared error is given by the conditional expectation
  $\fstar(x) \defn \Exs[Z \mid X = x]$.  Given a collection of samples
  $\{(x_i, z_i)\}_{i=1}^\numobs$, the goal of non-parametric
  regression is to produce an estimate $\fhat$ that is as close to
  $\fstar$ as possible.

\begin{figure}[h]
  \begin{center}
\begin{tabular}{ccc}
\raisebox{1.1in}{  \parbox{0.4\textwidth}{
    \begin{tabular}{cc}
      \widgraph{0.2\textwidth}{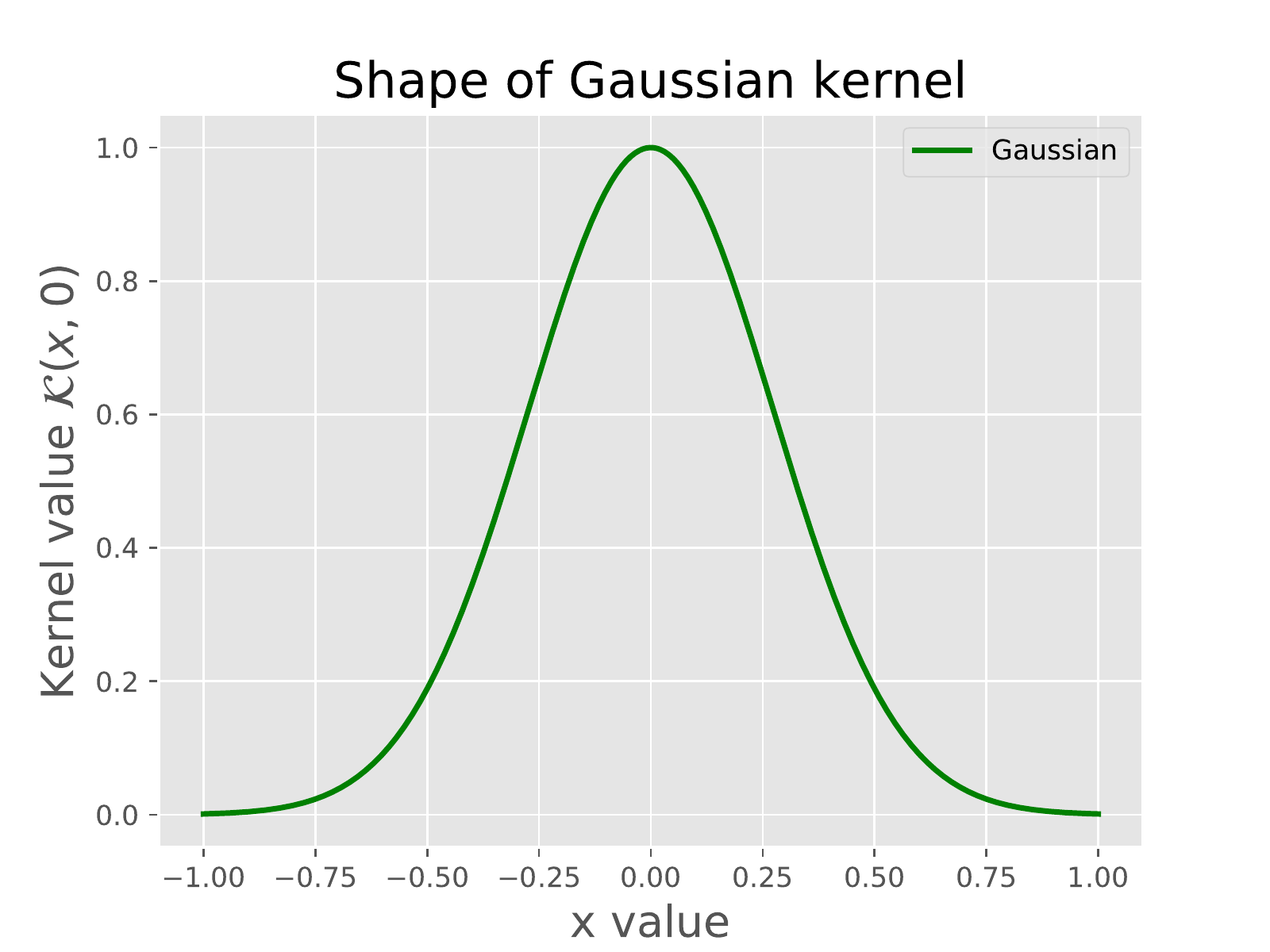} &
      \widgraph{0.2\textwidth}{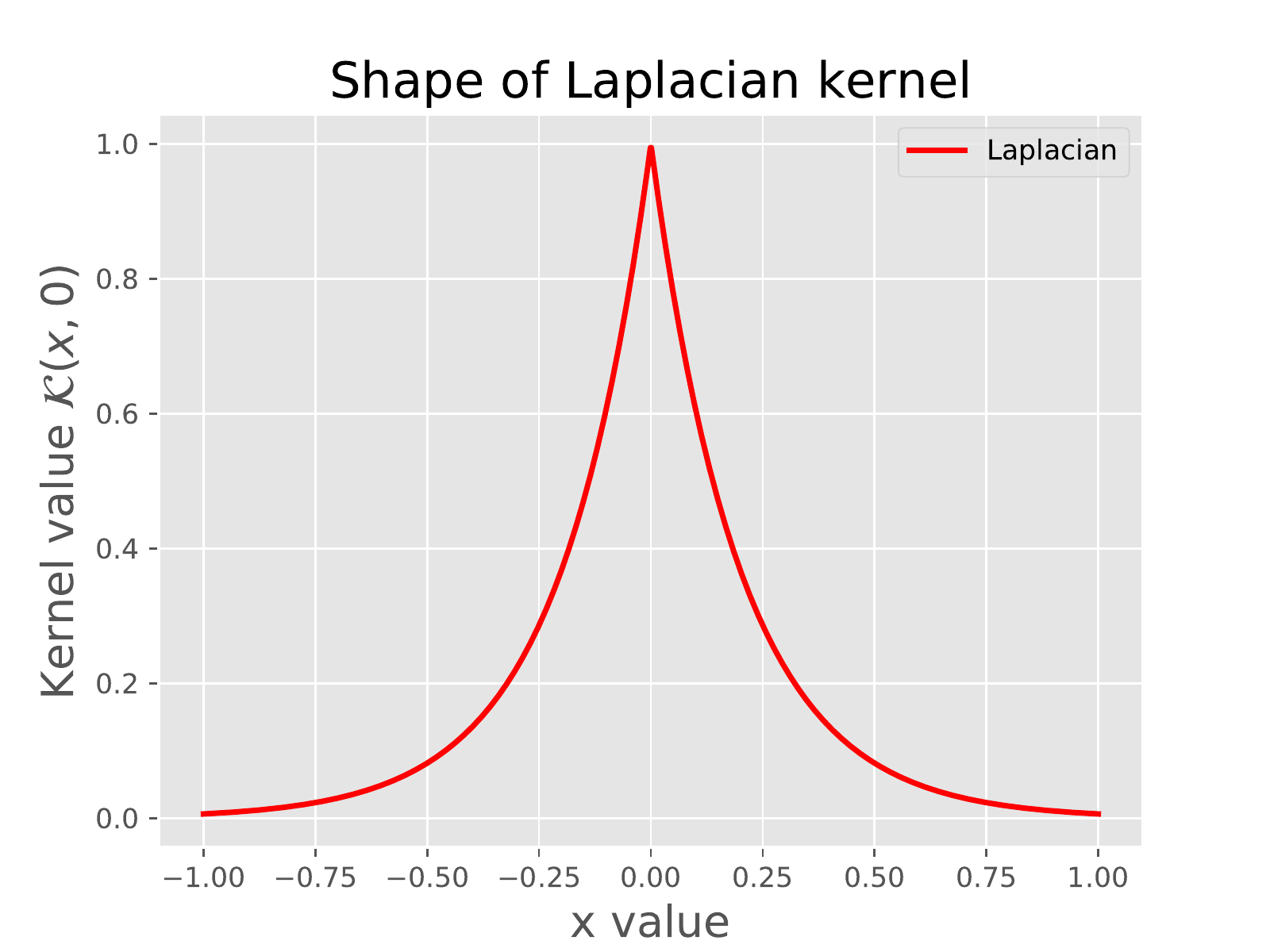} \\
      \widgraph{0.2\textwidth}{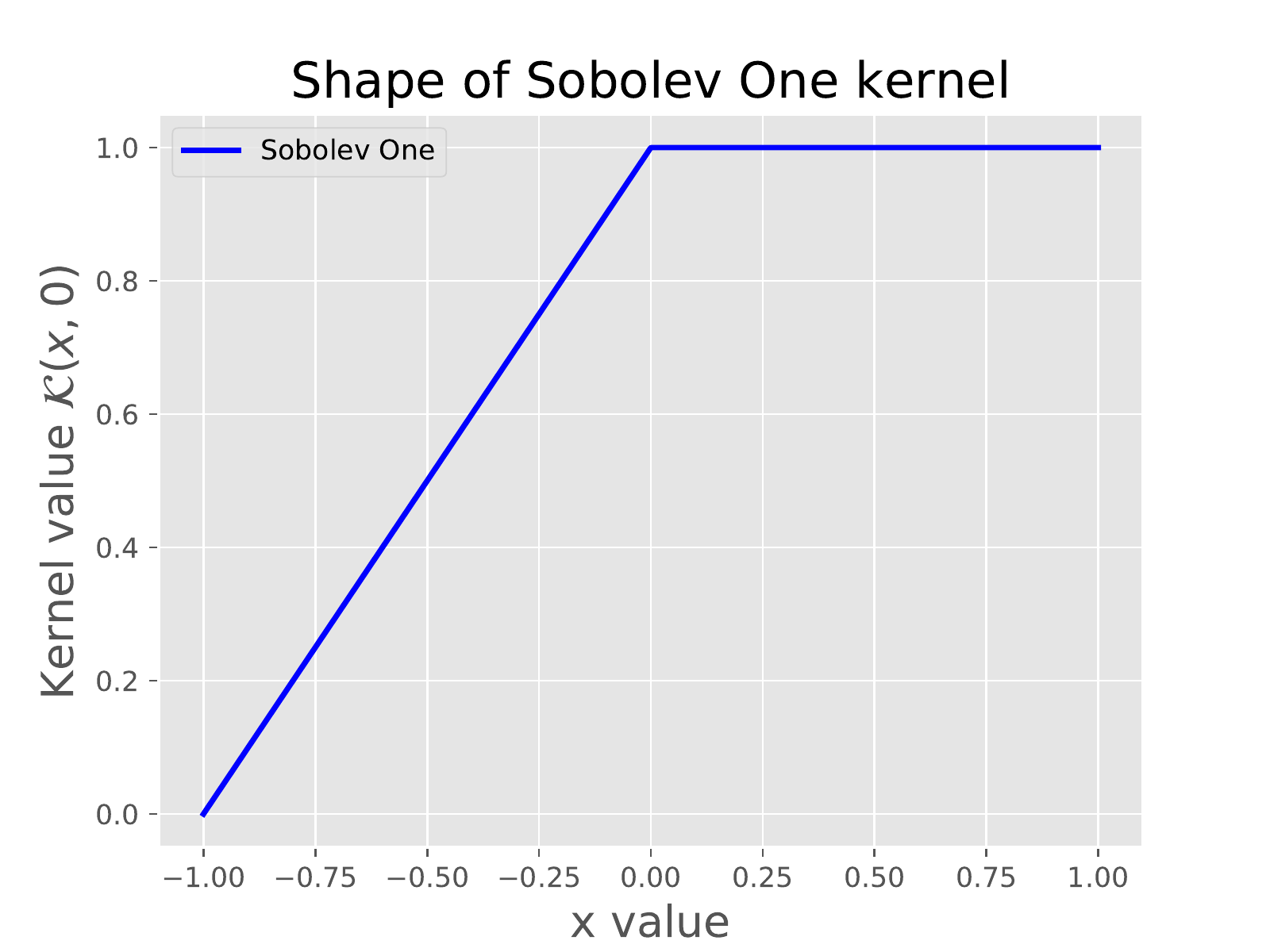} &
      \widgraph{0.2\textwidth}{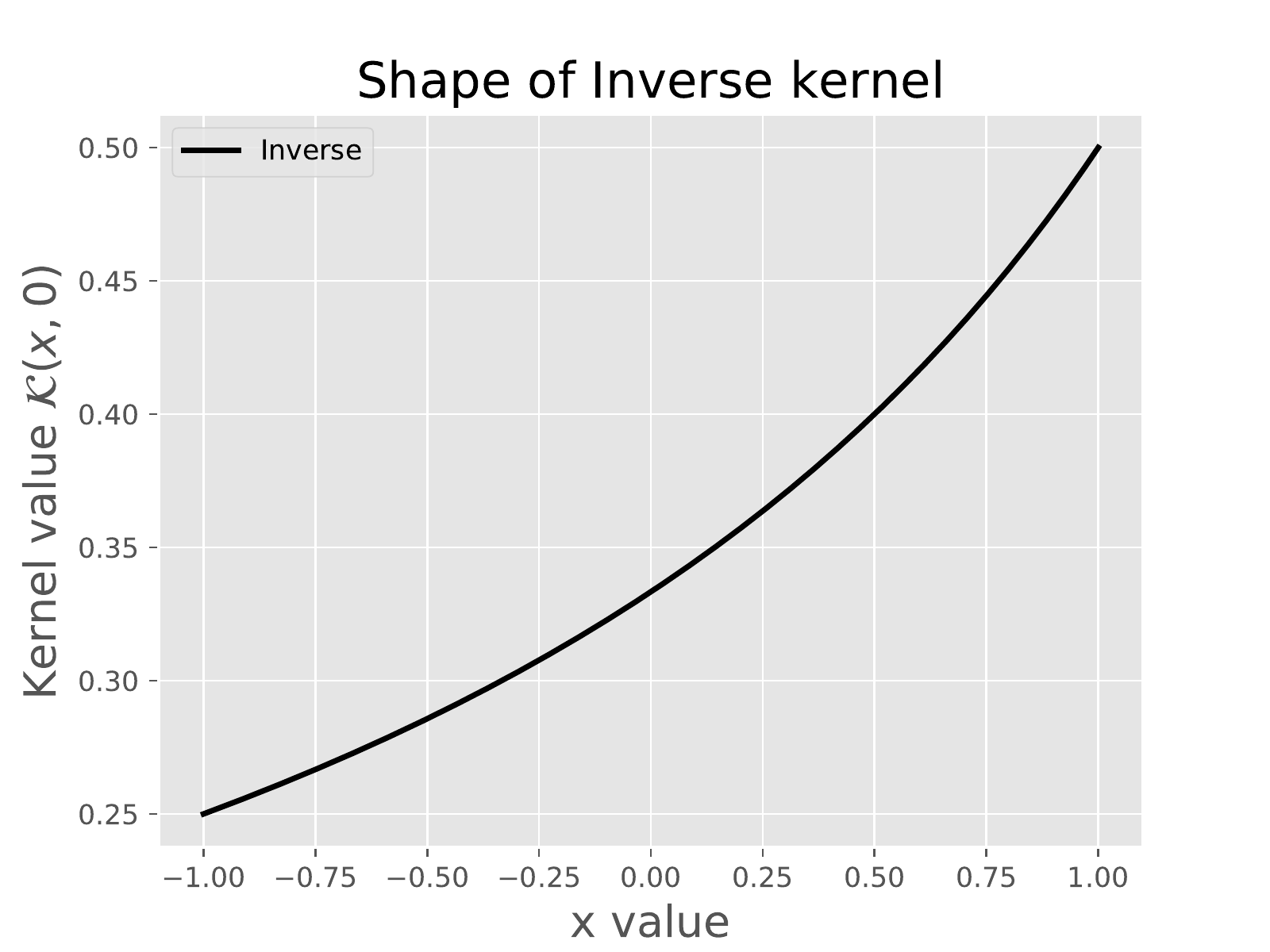}
    \end{tabular}
}}
& & \widgraph{0.5\textwidth}{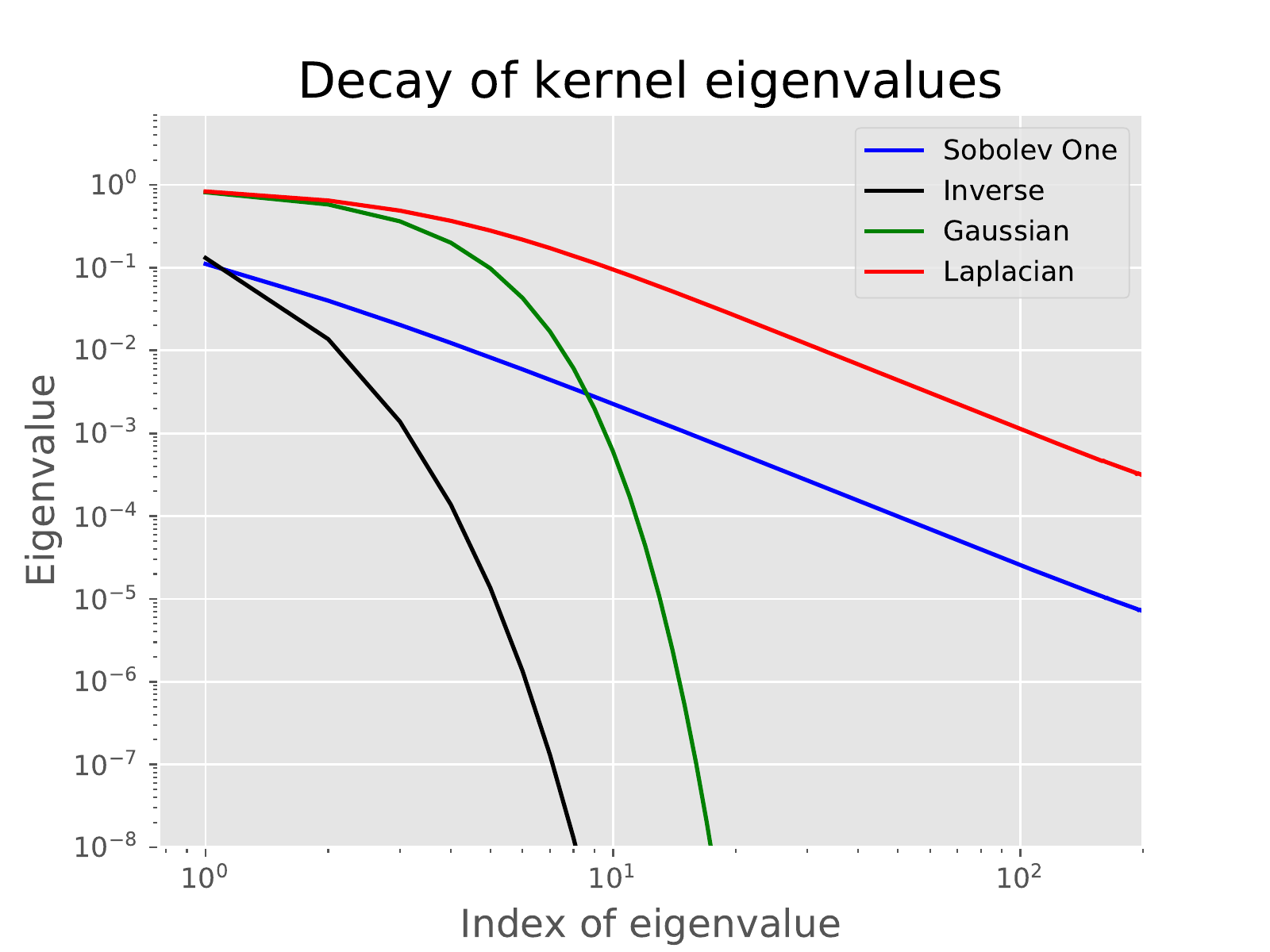} \\
(a) & & (b)
\end{tabular}
\caption{(a) Illustration of various kernel functions defined on
  $[-1,1] \times [-1, 1]$.  Each plot shows the kernel value
  $\PlainKer(x, 0)$ for $x \in [-1, 1]$.  (b) Illustration of the
  kernel eigenvalues $\{\mu_j\}_{j=1}^\numobs$ for kernel matrices $K$
  generated from the kernel functions in part (a).  Each log-log plot
  shows the eigenvalue versus the index: note how the Gaussian kernel
  eigenvalues decay at an exponential rate, whereas those of the
  Sobolev-One spline kernel decay at a polynomial rate.}
\label{FigKernels}
  \end{center}
\end{figure}

Assuming that the samples are i.i.d., we can rewrite our observations
in the form
\begin{align}
  \label{EqnOrigNonparaReg}
z_i = \fstar(x_i) + \gamma v_i, \qquad \mbox{for $i = 1, \ldots, \numobs$,}
\end{align}
where $v_i$ is an independent sequence of zero-mean noise variables
with unit variance.
A computationally attractive way of estimating $\fstar$ is to perform
least-squares regression over a \emph{reproducing kernel Hilbert
  space}, or RKHS for short~\cite{Aronszajn50,Kimeldorf71,Gu02,Wahba}.
Any such function class is defined by a symmetric, positive definite
kernel function $\KerPlain: \Xspace \times \Xspace \rightarrow \real$;
standard examples include the Gaussian kernel, Laplace kernel, and the
Sobolev (spline) kernels; see \autoref{FigKernels} for some
illustrative examples. Now suppose that $\fstar$ belongs to the RKHS
induced by the kernel $\KerPlain$, say with Hilbert norm
$\|\fstar\|_\Hil \leq R$.  In this case, the representer
theorem~\cite{Kimeldorf71} implies that the observation
model~\eqref{EqnOrigNonparaReg} is equivalent to
\begin{align*}
z & = \sqrt{n} K \alpha^* + \gamma v \qquad \mbox{for some $\alpha^* \in
  \real^\numobs$,}
\end{align*}
where $K \in \real^{\numobs \times \numobs}$ is the $\numobs \times
\numobs$ kernel matrix with entries $K_{ij} = \KerPlain(x_i, x_j)/n$
for each $i, j = 1, \ldots, \numobs$, and vector $v$ is a 
$\numobs$-dimensional vector formed by $v_i.$ The representer theorem 
and our choice of scaling ensures that $\|\fstar\|_\Hil^2 = 
(\alpha^*)^\top K \alpha^*$, meaning that $\alpha^*$ belongs to the 
ellipse of radius $R$ defined by the symmetric and PSD kernel matrix $K$.

Note that the matrix $K$ can be diagonalized as $K = U D U^\top$, where
$U$ is orthonormal, and \mbox{$D = \diag \{ \mu_1, \mu_2, \ldots,
  \mu_n\}$} is a diagonal matrix of non-negative eigenvalues.
Following this transformation, we arrive at an instance of the
standard ellipse model
\begin{align*}
  y & = \thetastar + w \qquad
  \mbox{where $w = \gamma U^\top v/\sqrt{\numobs}$, $\quad$ $y = U^\top
    z/\sqrt{\numobs}$,}
\end{align*}
and where $\thetastar = U^\top K \alpha^*$ belongs to the standard
ellipse~\eqref{EqnEllipse} defined by the eigenvalues of $K$.  Note
that the noise vector $w = \gamma U^\top v/\sqrt{\numobs}$ has zero-mean entries
each with standard deviation $\sigma = \gamma/\sqrt{\numobs}$.  The
entries of $w$ are not exactly Gaussian (unless the initial noise
vector $v$ was jointly Gaussian), but are often well-approximated by
Gaussian variables due to central limit behavior for large $n$.
\end{exas}


\subsection{Organization}

The remainder of this paper is organized as follows.  In
\autoref{SecSetup}, we introduce some background on
approximation-theoretic quantities, including the Gaussian width,
metric entropy, and the Kolmogorov width.  \autoref{SecMain} is
devoted to the statement of our main results, while
\autoref{SecEstimation} develops a number of their specific
consequences for ellipse estimation.  In \autoref{SecProofs}, we
provide the proofs of our main results, with more technical aspects of
the arguments provided in the appendices.


\section{Background}
\label{SecSetup}

Before proceeding to the statements of our main results, we introduce
some background on the notion of Gaussian width, Kolmogorov width, as
well as setting the estimation problem with ellipse constraint.

\subsection{Gaussian width}
\label{SecGaussianWidth}

Given a bounded subset $\Constraint \subset \real^\usedim$, the
\emph{Gaussian width} of $\Constraint$ is defined as
\begin{align}
  \label{EqnGaussWidth}
  \GWidth(\Constraint) \defn \Exs [\sup_{u\in \Constraint}
    \inprod{u}{\sgauss}] = \Exs \Big[\sup_{u\in \Constraint}
    \sum_{i=1}^\usedim \sgauss_i u_i \Big], \qquad \text{where }
  \sgauss_i \overset{\text{i.i.d.}}{\sim} \NORMAL(0, 1).
\end{align}
It measures the size of set $\Constraint$ in a certain sense.

It is also useful to define the classical notions of packing and
covering entropy.  An \emph{$\epsilon$-cover} of a set $\Constraint$
with respect to the $\ltwo{\cdot}$ metric is a discrete set
$\{\theta^1,\ldots,\theta^N\} \subset \Constraint$ such that for each
$\theta \in \Constraint$, there exists some $i\in \{1,\ldots, N\}$
satisfying $\ltwo{\theta - \genparam^i} \leq \epsilon$.  The
\emph{$\epsilon$-covering number} $\Ncover(\epsilon, ~\Constraint)$ is
the cardinality of the smallest $\epsilon$-cover, and the logarithm of
this number $\log \Ncover(\epsilon, ~\Constraint)$ is called the
\emph{covering metric entropy} of set $\Constraint$.

Similarly, an \emph{$\epsilon$-packing} of a set $\Constraint$ is a set
$\{\genparam^1, \ldots, \genparam^M\} \subset \Constraint$ satisfying
$\ltwo{\genparam^i - \genparam^j} > \epsilon$ for all $i \ne j$.  The
size of the largest such packing is called the
\emph{$\epsilon$-packing number} of $\Constraint$, which we denote by
$\Mpack(\epsilon, \Constraint)$.  It is related to the (covering)
metric entropy by the inequalities
\begin{align*}
\log \Mpack(2\epsilon, \Constraint) \, \le \, \log \Ncover(\epsilon,
\Constraint) \, \le \, \log \Mpack(\epsilon, \Constraint).
\end{align*}
For this reason, we use the term metric entropy to refer to either the
covering or packing metric entropy, since they differ only in constant
terms.

The connection between Gaussian width and metric entropy is
well-studied (e.g. \cite{Dudley67,talagrand2014upper,Wai17}).
For our future discussion, we collect a few results here as reference.   
First, Dudley's entropy integral~\cite{Dudley67}
is an upper bound for the Gaussian width---viz.
\begin{align}
\GWidth(\Constraint)
\,\le \, c \int_0^{\diam(\Constraint)}
\sqrt{\log \Ncover(\epsilon, \Constraint)} \mathop{d\epsilon},
\end{align}
for some universal constant $c>0$.
This upper bound also holds for more general sub-Gaussian processes.
Dudley's bound can be much looser than the more refined bounds
obtained through Talagrand's generic chaining, which are tight up to a
universal constant~\cite[Thm. 2.4.1]{talagrand2014upper}.  For
Gaussian processes like ours, Sudakov minoration (e.g.,
~\cite[Thm. 13.4]{BouLugMas13}) provides a lower bound on the Gaussian
width.
\begin{align}
\label{EqnSudakovMinoration}
\GWidth(\Constraint)
\, \ge \, \sup_{\epsilon > 0}\, c \, \epsilon \, \sqrt{\log \Mpack(\epsilon, \Constraint)}.
\end{align}
Although we do not directly use this lower bound when proving our main
lower bound (\autoref{ThmLowerBound}) below,
we follow its spirit by constructing a large collection of well-separated points.


\subsection{Kolmogorov width}
In this section, we define the Kolmogorov width and briefly review its
properties.  This geometric quantity plays the central role in our
main results.

For a given compact set $\Constraint \subset \real^\usedim$ and
integer $k \in [\usedim]$, the \emph{Kolmogorov $k$-width} of
$\Constraint$ is given by
\begin{align}
  \label{EqnKolmogorov}
\kwidth_k (\Constraint) \defn \min_{\trunc{k} \in \mathcal{P}_k}
\max_{\theta \in \Constraint} \|\theta - \trunc{k} \theta\|_2,
\end{align}
where $\LPset{k}$ denotes the set of all $k$-dimensional orthogonal
linear projections, and $\trunc{k} \theta$ denotes the projection of 
$\theta$ to the corresponding $k$-dimensional linear space. 
Any projection $\trunc{k}$ achieving the minimum
in expression~\eqref{EqnKolmogorov} is said to be an \emph{optimal
  projection} for $\kwidth_k(\Constraint)$.  Note that the Kolmogorov
width $\kwidth_k(\Constraint)$ is a non-increasing function of $k$,
meaning that
\begin{align*}
\max \limits_{\theta \in \Constraint} \ltwo{\genparam} =
\kwidth_0(\Constraint) \geq \kwidth_1(\Constraint) \geq \ldots \geq
\kwidth_d(\Constraint) = 0.
\end{align*}
We refer the readers to the book by \citet{pinkus2012n} for more
details on the Kolmogorov width and its properties.



\section{Main results}
\label{SecMain}
Let us first define the notion of localized Gaussian width formally,
and then turn to the statement of our main results.

\subsection{Localized Gaussian width}

\begin{figure}[ht]
	\centering
	\includegraphics[width=0.5\textwidth]{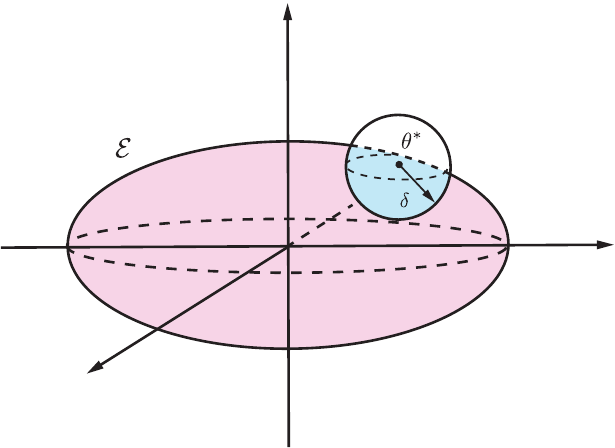}
	\caption{An illustration of the set $\Ellipse_{\at} \cap
          \Ball(\delta)$.  It is the intersection of the ellipse with
          Euclidean ball centered at $\at$, and thus varies according
          to the local geometry of the ellipse.}
	\label{fig:ellips}
\end{figure}

Let $\Ball(\delta)$ denote the Euclidean ball of radius $\delta$, and
for a given vector $\thetastar \in \Ellipse$, define the shifted
ellipse $\Ellipse_\at \defn \big \{ \theta - \at \mid \theta \in
\Ellipse \big \}$.  The \emph{localized Gaussian width at $\at$ and
  scale $\delta$} is defined as
\begin{align}
\label{EqnGaussWidthLocBall}
\GWidth(\Ellipse_{\at} \cap \, \Ball(\delta))
= \Exs \sup_{\Delta \in \Ellipse_{\at} \cap \Ball(\delta)}
\inprod{\sgauss}{\Delta}.
\end{align}
Note that this quantity is simply the ordinary Gaussian width of the
set $\Ellipse_{\at} \cap \Ball(\delta)$, and we say that it is
localized since the Euclidean ball restricts it to a neighborhood of
$\at$.  See \autoref{fig:ellips} for an illustration of this set.

We note that localized forms of Gaussian and Rademacher complexity are
standard in the literature on empirical processes
(e.g.,~\cite{Bar05,Kolt06}), where it is known that they are needed to
obtain sharp rates.  In the case of least-squares estimation over
convex sets, there is an extremely explicit connection between the
localized Gaussian width and the associated estimation
error~\cite{vandeGeer,chatterjee2014new,Wai17}; we describe this
relationship in more detail in \autoref{SecEstimation} and
\autoref{SecEstimationProof}.

Our main results, to be stated in the following subsections, provide
conditions under which we can provide a sharp characterization of the
localized Gaussian width~\eqref{EqnGaussWidthLocBall} in terms of the
Kolmogorov width.

\subsection{Upper bound on the localized Gaussian width}
\label{SecUB}

In order to state our first main result, we introduce an
approximation-theoretic quantity having to do with the quality of a
given $k$-dimensional projection. For a given integer $k \in \{1,
\ldots, \usedim \}$ and any $k$-dimensional linear projection
$\Proj_{k}$, let us define the set
\begin{align}
\Gammaclass(\genparamstar, \delta, \Proj_k) \defn \braces*{ \gamma \in
  \real^\usedim \;\middle|\; \gamma > 0, \; \sup_{\Delta \in
    \Ellipse_{\at} \cap \Ball(\delta)} \sum_{i=1}^{\usedim}
  \frac{[\Delta_i - (\Proj_{k} \Delta)_i]^2}{\gamma_i} \le 1 }.
\end{align}
Here $\gamma > 0$ means that $\gamma_i > 0$ for each coordinate 
$i=1,\ldots,\numobs.$
It can be verified that the set $\Gammaclass(\genparamstar, \delta,
\Proj_k)$ is always non-empty since the constant vector \mbox{$\gamma
  = \sqrt{\mu_1}\delta \, \mathbf{1}$} always belongs to it.  (Here
$\mathbf{1}$ denotes the vector of all ones.)  To provide some
intuition for this definition, the vector $\Delta - \Proj_k(\Delta)$
corresponds to the error incurred by using the subspace associated
with $\Proj_k$ to approximate $\Delta$.  The positive vector $\gamma \in
\real^\usedim$ allows us to weight the entries of this error vector in
computing the Euclidean norm of the weighted error.

\noindent We are now ready to state an upper bound on the localized
Gaussian width.
\begin{theos}
\label{ThmUpperBound}
Given any $\delta > 0$, projection tuple $(k,\Proj_k)$, and vector
$\genparamstar \in \Ellipse$, we have
\begin{align}
\label{EqnUB}
\GWidth(\Ellipse_{\at} \cap \Ball(\delta)) \leq \delta \, \sqrt{k} +
\inf_{\gamma \in \Gammaclass(\genparamstar, \delta, \Proj_k)} \;
\sqrt{\sum_{i=1}^\usedim \gamma_i}.
\end{align}
\end{theos}
\noindent See \autoref{SecPfUB} for the proof of this result.

Note that \autoref{ThmUpperBound} holds for any dimension and
projection pair $(k,\Proj_k)$.  Often the case, we can choose a
specific pair for which the set $\Gammaclass(\genparamstar, \delta,
\Proj_k)$ is easy to characterize. In particular, given any fixed
$\delta > 0$, let us define the \emph{critical dimension}
\begin{align}
\label{Eqnkcrit}
\lwk(\genparamstar, \delta) \defn \arg \min_{k = 1, \ldots, \usedim}
\Big \{ \kwidth_{k}\Big(\Ellipse_{\at} \cap \Ball\big((1 -
\tinyconst)\delta\big)\Big) \leq \frac{9}{10} \delta \Big \},
\end{align}
for some constant $\tinyconst \in (0,0.1)$. In words, this integer is
the minimal dimension for which there exist a $\lwk$-dimensional
projection that approximates a neighborhood of the re-centered ellipse
to $\frac{9}{10} \delta$-accuracy.\footnote{The constants $\tinyconst$
  and $9/10$ are chosen for the sake of convenience in the proof, but
  other choices of these quantities (which both must be strictly less
  than $1$) are also possible.}  Although our notation does not
explicitly reflect it, note that $\lwk(\genparamstar, \delta)$ also
depends on the \mbox{ellipse $\Ellipse$.}

Given the integer $\upk \equiv \upk(\thetastar, \delta)$, we let
$\Proj_{\upk} \in \ProjClass{\upk}$ denote the minimizing projection
in the definition~\eqref{EqnKolmogorov} of the width, and note that
for any vector $\Delta$, the error associated with this projection is
given by $\Delta - \Proj_{\upk}(\Delta)$.  It can be seen in our later
examples, this particular choice $(\upk, ~\Proj_{\upk})$ often yields
tight control of the localized Gaussian width.  So as to streamline
notation, we adopt $\Gammaclass(\genparamstar, \delta)$ as a short
hand for $\Gammaclass(\genparamstar, \delta, \Proj_{\upk})$.

\paragraph{Regularity assumption:}
For many ellipses encountered in practice, the first term in the upper
bound~\eqref{EqnUB} dominates the second term involving the set
$\Gammaclass$.  In order to capture this condition, we say the ellipse
$\Ellipse$ is \emph{regular at $\genparamstar$} if there exists some
pair $(k,\Proj_k)$ such that
\begin{align}
\label{EqnRegularity}
\inf_{\gamma \in \Gammaclass(\genparamstar, \delta, \upk,
  \Proj_{\upk})} \sum_{i=1}^\usedim \gamma_i \,\leq \, c\; \delta^2 \,
k \qquad \mbox{for all $\delta > 0$.}
\end{align}
Here $c < \infty$ is any universal constant.  When this condition
holds, \autoref{ThmUpperBound} implies the existence of another
universal constant $c'$ such that
\begin{align}
\GWidth(\Ellipse_{\at} \cap \Ball(\delta)) \leq c' \; \delta \,
\sqrt{k} \qquad \mbox{for all $\delta > 0$.}
\end{align}
As is shown in \autoref{SecKernelRegularity}, the regularity
condition~\eqref{EqnRegularity} is a generalization of a condition
previously introduced by~\citet{yang2015randomized} in the context of
kernel ridge regression, and it holds for many examples encountered in
practice.

\noindent As a direct consequence of \autoref{ThmUpperBound}, the
following corollary holds.
\begin{cors}
If the regularity assumption~\eqref{EqnRegularity} is satisfied with dimension and
projection pair $(\upk, \Proj_{\upk})$, then the localized Gaussian width
satisfies
\begin{align}
\label{EqnUBsimple}
\GWidth(\Ellipse_{\at} \cap \Ball(\delta)) \leq \constup \; \delta \,
\sqrt{\upk} \qquad \mbox{for all $\delta > 0$.}
\end{align}
\end{cors}

\noindent Let us illustrate the regularity
condition~\eqref{EqnRegularity} and associated consequences of
\autoref{ThmUpperBound} with some examples.

\begin{exas}[Gaussian width of the Euclidean ball]
\label{ExampleBall}
We begin with a simple example: suppose that the ellipse $\Ellipse$ is
the Euclidean ball in $\real^\usedim$, specified by the aspect ratios
$\mu_j = 1$ for all $j = 1, \ldots, \usedim$, and let us use
\autoref{ThmUpperBound} to upper bound the Gaussian width at
$\thetastar = 0$.  For $\delta \in (0,\frac{1}{1 - \tinyconst})$ and
any integer $k < \usedim$, we have \mbox{$\kwidth_{k}(\Ellipse_{\at}
  \cap \Ball((1 - \tinyconst) \delta)) = (1 - \tinyconst) \delta$,}
because any $k$-dimensional projection must neglect at least one
coordinate.  Since $1 - \tinyconst > 9/10$, we conclude that $\lwk(0,
\delta) = \usedim$ for all $\delta \in (0,\frac{1}{1 - \tinyconst})$.
With this choice of $\lwk$, there is no error in the projection,
meaning that $\inf_{\gamma \in \Gammaclass(\genparamstar, \delta)}
\sum_{i=1}^\usedim \gamma_i = 0$.  Consequently, the regularity
condition~\eqref{EqnRegularity} certainly holds, so that
\autoref{ThmUpperBound} implies that
\begin{align*}
  \GWidth(\Ellipse_{\at} \cap \Ball(\delta)) \leq c' \delta
  \sqrt{\usedim}.
\end{align*}
 In fact, a direct calculation yields that
$\GWidth(\Ellipse_{\at} \cap \Ball(\delta)) = \delta (\sqrt{\usedim} -
o(1))$, where $o(1)$ is a quantity tending to zero as $\usedim$ grows
(e.g.,~\cite{Wai17}).  Consequently, our bound is asymptotically sharp
up to the constant pre-factor in this special case.
\end{exas}

\vspace*{.2in}

\noindent We now turn to a second example that arises in non-parametric
regression and density estimation under smoothness constraints: \\

\begin{exas}[Gaussian width for Sobolev ellipses]
\label{ExampleSobolev}
Now consider an ellipse $\Ellipse$ defined by the aspect ratios
$\eig{j} = c j^{-2 \smoothpar}$, where $\smoothpar > 1/2$ is a
parameter.  Ellipses of this form arise when studying non-parametric
estimation problems involving functions that are $\smoothpar$-times
differentiable with Lebesgue-integrable
$\smoothpar$-derivative~\cite{Tsybakov09}.  Let us again use
\autoref{ThmUpperBound} to upper bound the localized Gaussian width at
$\thetastar = 0$.  From classical results on Kolmogorov widths of
ellipses~\cite{pinkus2012n} (see also \cite[Sec. 4.3]{wei2017testing}),
we know that $\kwidth_{k}(\Ellipse_{0}) =
\sqrt{\mu_{k+1}}$.  Taking into account the intersection with the
Euclidean ball, we find that
\begin{align}
\label{EqnSobolevArgument}
  \kwidth_{k}(\Ellipse_{\at} \cap \Ball((1 - \tinyconst) \delta) =
  \min \big \{ \sqrt{\mu_{k+1}}, (1 - \tinyconst) \delta \big \},
\end{align}
valid for any $\delta \in (0, \frac{1}{1 - \tinyconst}\sqrt{\mu_1})$.
Since $1 - \tinyconst > 9/10$, we conclude that
\begin{align*}
  \lwk(0, \delta) = \arg\min \{\sqrt{\mu_{k+1}} \leq \frac{9}{10} \delta \}
 = \lceil (\frac{10\sqrt{c}}{9\delta})^{1/\smoothpar} \rceil,
\end{align*}
again valid for all $\delta \in (0,\frac{1}{1 -
  \tinyconst}\sqrt{\mu_1})$.  
Here the last inequality uses the fact that $\eig{j} = c j^{-2 \smoothpar}$.

This argument also shows that the
corresponding projection subspace is spanned by the first $\lwk$
standard orthogonal vectors $\{e_i\}_{i=1}^{\lwk}$.  With this
projection, any feasible vector $\gamma \in \Gammaclass(\genparamstar,
\delta)$ satisfies $\gamma_i \geq \mu_{i} \indic\{i > \lwk(0,
\delta)\}$, meaning that
\begin{align}
\label{EqnRegAtZero}
  \inf_{\gamma \in \Gammaclass(\genparamstar, \delta)} \sum_{i=1}^\usedim \gamma_i
= \sum_{j = \lwk+1}^\usedim \mu_{j}
= c \sum_{j = \lwk+1}^\usedim j^{-2\smoothpar} \leq  c\int_{\lwk+1}^{\infty} t^{-2\smoothpar} \mathop{dt}
= c \delta^{2-1/\smoothpar}.
\end{align}
On the other hand, we also have $\delta^2 \lwk(0, \delta) \asymp \delta^{2
  - 1/\smoothpar}$, so there exists
some constant $c'$, such that $\inf_{\gamma \in
  \Gammaclass(\genparamstar, \delta)} \sum_{i=1}^\usedim \gamma_i \leq
c' \delta^2 \lwk(0, \delta)$ which validates the regularity
condition~\eqref{EqnRegularity}.  Therefore,
\autoref{ThmUpperBound} guarantees that
\begin{align}
\label{EqnUBSobolev}
\GWidth(\Ellipse_{0} \cap \Ball(\delta)) \leq c'' \; \delta^{1-
  (1/2\smoothpar)}.
\end{align}
In fact, the above bound~\eqref{EqnUBSobolev} can be shown to be tight
up to a constant pre-factor.  See the discussion following
\autoref{CorMetEnt} in the sequel for further details.
\end{exas}


\subsection{Lower bound on the localized Gaussian width}
\label{SecLB}

Thus far, we have derived an upper bound for the localized Gaussian
width. In this section, we use information-theoretic methods to prove
an analogous lower bound on the localized Gaussian width.  This lower
bound involves both the critical dimension $\kind(\genparamstar,
\delta)$, as previously defined in equation~\eqref{Eqnkcrit}, and also
a second quantity, one which measures the proximity of $\thetastar$ to
the boundary of the ellipse. More precisely, for a given $\thetastar
\in \Ellipse$, define the mapping $\Rfun: \real_+ \rightarrow \real_+$
via
\begin{align}
\label{EqnRfun}
\Rfun(\delta) & = \begin{cases} 1 & \mbox{if $\delta >
    \|\thetastar\|_2 / (1 - \tinyconst)$} \\
    1 \wedge \min \Big \{ r \geq 0 \, \mid
  \delta^2 \leq \frac{1}{(1 - \tinyconst)^2}
  \sum_{i=1}^\usedim \frac{r^2}{(r+\mu_i)^2} (\at_i)^2
  \Big \} & \mbox{otherwise.}
  \end{cases}
\end{align}
As shown by \citet{wei2017testing}, this mapping is
well-defined, and has the limiting behavior $\Rfun(\delta) \rightarrow
0$ as $\delta \rightarrow 0^+$; for completeness, we include the verification
of these claims in \autoref{AppRfun}, along with a sketch of the function.
Let us denote
$\Rfun^{-1}(x)$ as the largest positive value of $\delta$ such that
$\Rfun(\delta) \leq x$. Note that by this definition, we have
$\InvRfun(1) = \infty$.

Recall that the elliptical norm on $\real^\usedim$ is defined via 
$\enorm{\genparam}^2 \defn \sum_{j=1}^\usedim \frac{\genparam_j^2}{\eig{j}}$.
We are now ready to state our lower bound for the localized Gaussian
width.
\begin{theos}
\label{ThmLowerBound}
There exist universal constants $\constlw, c > 0$ such that for all
$\thetastar \in \Ellipse$
\begin{align}
\GWidth(\Ellipse_{\at} \cap \Ball(\delta)) \geq \constlw \; \delta \, \sqrt{1
  - \enorm{\genparamstar}^2} \, \sqrt{\lwk(\genparamstar, \delta)},
\qquad \mbox{for all $\delta \in \Big(0, c \Rfun^{-1} \big(
  (\enorm{\thetastar}^{-1} - 1)^2 \big) \wedge \sqrt{\mu_1} \Big)$.}
\end{align}
\end{theos}
\noindent See \autoref{SecPfLB} for the proof of this theorem.

We remark that the regularity condition~\eqref{EqnRegularity} is not
necessary for this result to hold.  Besides, in order to understand the
inequality \mbox{$\delta < c \Rfun^{-1} \big((\enorm{\thetastar}^{-1}
  - 1)^2 \big)$}, it is equivalent to ask for \mbox{$\enorm{\at} <
\frac{1}{1 + \sqrt{\Rfun(\delta/c)}}$}. We assume this since it is not our 
primary interest to study the case when $\at$ is sufficiently close to the
boundary of the ellipse.  Concretely, if we assume that $\enorm{\at}
\leq 1/2$, then $(\enorm{\thetastar}^{-1} - 1)^2 \geq 1$ therefore
$\InvRfun \big ( (\enorm{\thetastar}^{-1} - 1)^2 \big) = \infty$.


\subsection{Some consequences}

One useful consequence of \autoref{ThmUpperBound} and
\autoref{ThmLowerBound} is in providing sufficient conditions for
tight control of the localized Gaussian width.  If the ellipse
$\Ellipse$ is regular at $\genparamstar$, then the above theorems
imply the localized Gaussian width~\eqref{EqnGaussWidthLocBall} is
equivalent to $\delta \sqrt{\upk(\genparamstar, \delta)}$ up to a
multiplicative constant. Specifically, we have the sandwich relation
\begin{align}
\label{EqnGWidthSandwich}
\constlw \delta \sqrt{\lwk(\genparamstar, \delta)} \, \le \,
\GWidth(\Ellipse_{\at} \cap \Ball(\delta))  \,\le \, \constup \delta
\sqrt{\upk(\genparamstar, \delta)},
\end{align}
for some positive constants $\constlw$ and $\constup$.

Recall our earlier calculation from \autoref{ExampleBall}, where we
showed that the localized Gaussian width scales as $\delta
\sqrt{\usedim}$, up to multiplicative constants.  The sandwich
relation~\eqref{EqnGWidthSandwich} shows that this same scaling holds
more generally with $\usedim$ replaced by $\upk(\genparamstar,
\delta)$.  Thus, we can think of $\upk(\genparamstar, \delta)$
corresponding to the ``effective dimension'' of the set $\Ellipse_\at
\cap \Ball(\delta)$.

It is worthwhile pointing out that our results have a number of
corollaries, in particular in terms of how local Gaussian widths and
Kolmogorov widths are related to metric entropy.  Recall the notion of
the metric (packing) entropy $\log M$ as previously defined
\autoref{SecGaussianWidth}.  The following corollary provides a
sandwich for $\lwk(\genparamstar, \delta)$ in terms of the metric
entropy of the set $\Ellipse_\at \cap \Ball(\delta)$.
\begin{cors}
\label{CorMetEnt}
There are universal constants $c_j > 0$ such that for any pair $(\at,
\Ellipse)$ satisfying the regularity condition~\eqref{EqnRegularity},
we have
\begin{align}
\label{EqnMetEnt}
  c_1 \log \Mpack \big(\frac{\delta}{2}, \Ellipse_{\at} \cap
  \Ball(\delta) \big) \; \overset{(i)}{\leq} \; \kind(\genparamstar,
  \delta) & \; \overset{(ii)}{\leq} \; c_2 \, \log \Mpack \big(c_0
  \delta, \Ellipse_{\at} \cap \Ball(\delta) \big) \quad \mbox{for all
    $\delta \in \big(0, \frac{1}{e} \big)$. $\quad$}
\end{align}
\end{cors}
\noindent See~\autoref{AppProofCorMetEnt} for the proof.  The lower
bound (i) is a relatively straightforward consequence of Sudakov's
inequality~\eqref{EqnSudakovMinoration}, when combined with our
results connecting the Kolmogorov and Gaussian widths.  The upper
bound (ii) requires a lengthier argument.

Recall that in \autoref{ExampleSobolev}, we argued that for the
Sobolev ellipse with smoothness $\alpha > 1/2$, the Kolmogorov width
at $\thetastar = 0$ is given by $\upk(0, \delta) = c \,
(1/\delta)^{(1/\alpha)}$.  Combining this calculation with
\autoref{CorMetEnt}, we find that $\log \Mpack \big(\delta/2,
\Ellipse_{\at} \cap \Ball(\delta) \big) = (1/\delta)^{1/\alpha}$
up to a multiplicative constant.
This is a known fact that can be verified by constructing explicit
packings of these function classes, but it serves to illustrate the
sharpness of our results in this particular context.


\section{Consequences for estimation}
\label{SecEstimation}

In the previous section, we established upper and lower bounds on the
localized Gaussian width in \autoref{ThmUpperBound} and
\autoref{ThmLowerBound}.  We now turn to some consequences of these
bounds, in particular for the problem of constrained least-squares
estimation.

In particular, suppose we are given observations
$y \sim \NORMAL(\genparamstar, \noisestd^2 \IdMat_{n})$
with $\genparam \in \Ellipse$ according
to the earlier model~\eqref{EqnGaussianModel},
and we consider the constrained least squares estimator (LSE)
\begin{align}
\label{EqnLS}
\genparamhat \defn \arg\min_{\genparam \in \Ellipse} \ltwo{y -
  \genparam}^2.
\end{align}
Let us assume that the ellipse $\Ellipse$ is regular at
$\genparamstar$, so that the localized Gaussian width satisfies the
bounds~\eqref{EqnGWidthSandwich} with constants $\constlw$ and
$\constup$.  Connecting the error $\ltwo{\genparamhat -
  \genparamstar}$ to these Gaussian width bounds involves the
following two functions
\begin{align}
\label{EqnCritBounds}
\CritFunUp(\delta)  \defn \frac{\delta^2}{2} - \noisestd \constlw \delta
\sqrt{\lwk(\at, \delta)}, \quad \mbox{and} \quad
\CritFunLw(\delta)  \defn \frac{\delta^2}{2} - \noisestd \constup
\delta \sqrt{\upk(\at, \delta)},
\end{align}
with the critical dimension defined in expression~\eqref{Eqnkcrit}.

Let us consider the fixed point equation
\begin{align}
\label{EqnFixedPt}
\delta &= \constlw \noisestd \sqrt{\lwk(\at, \delta)}
\qquad
\text{for }
\delta \le c \Rfun^{-1} \big(
  (\enorm{\thetastar}^{-1} - 1)^2 \big) \wedge \sqrt{\mu_1}.
\end{align}
Since $\delta \mapsto \lwk(\delta)$ is a non-increasing function of
$\delta$ (see \citet[Appendix D.1]{wei2017testing}) while $\delta
\mapsto \delta$ is increasing, if this fixed point
problem~\eqref{EqnFixedPt} has a solution, then the solution is unique
and we denote it as $\delcrit$.

We can now give a precise statement relating the estimation rate of
$\genparamhat$ to the solution $\delcrit$ of the fixed point
equation~\eqref{EqnFixedPt}.



\begin{props}[Least squares on ellipses]
\label{PropEllipseEstimation}
Let $\Ellipse$ be regular at $\genparamstar$, and let
$\delcrit$ be the solution to the fixed point
problem~\eqref{EqnFixedPt}.
Suppose furthermore the following conditions hold
\begin{enumerate}[(a)]
  \item The function $\CritFunLw$ is unimodal in $\delta$.

  \item There exists a constant $c_1 \in (0,1)$ such that
    $\constup^2\upk(\delta) \leq \frac{1}{4 c^2_1} \constlw^2\lwk(\delcrit)$ for
    $\delta = c_1 \delcrit$,

  \item There exists a constant $c_2 > 1$ such that $\delta \geq 2
    \noisestd \constup \sqrt{\upk(\delta)}$ for $\delta = c_2 \delcrit$.
\end{enumerate}
Then the error of the least squares estimator~\eqref{EqnLS} satisfies
\begin{align}
  \label{EqnHighProb}
c \delcrit \leq \ltwo{\genparamhat - \genparamstar} \leq c' \delcrit,
\qquad \text{with prob. } \ge 1 - 3 \exp(-c'' \delcrit^2 /
\noisestd^2),
\end{align}
for some constants that depend only on $c_1$ and $c_2$.
\end{props}
\noindent See \autoref{SecEstimationProof} for the proof of this
result.

Note that this result is stated for the ellipse $\Ellipse(\bigrad)$ 
with $\bigrad=1$.  For arbitrary $\bigrad$ one can easily rescale to 
obtain similar results; see equation~\eqref{EqnRescaling} in 
\autoref{SecReduction} for more detail.
When we say $\CritFunLw$ is unimodal, we mean that there is some $t$
such that $\CritFunLw$ is nondecreasing for $\delta < t$ and
nonincreasing for $\delta > t$.

Equation~\eqref{EqnHighProb} provides a high probability bound on the
least-squares error.  If furthermore $\delcrit \gtrsim \noisestd$,
then we are also guaranteed that the mean-squared error is sandwiched
as
\begin{align}
\label{EqnExpectationBound}
c \delcrit^2 \; \leq \; \Exs \ltwo{\genparamhat - \genparamstar}^2 \;
\leq \; c' \delcrit^2
\end{align}
for some universal constants $(c, c')$.

We claim the conditions of \autoref{PropEllipseEstimation} are
relatively mild.  Note that the related function $g(t) \defn
\frac{\delta^2}{2} - \noisestd \GWidth(\Ellipse_{\at} \cap
\Ball(\delta))$ is strongly convex \cite[Thm. 1.1]{chatterjee2014new},
as mentioned in Appendix~\ref{SecReduction}.  So it is reasonable to
believe that its approximation $\CritFunLw$ is unimodal.  Moreover,
the assumptions (b) and (c) essentially assert that $\CritFunLw$ does
not change too drastically at two points $c_1 \delcrit$ and $c_2
\delcrit$ close to the critical radius $\delcrit$.
In the next section, we will check these assumptions for different
examples.

Note that fixed point problem~\eqref{EqnFixedPt}
can be viewed as a kind of a \emph{critical equation}
(e.g., \cite[Ch. 13]{Wai17} and \cite{yang2015randomized}),
whose solution $\delcrit$ we call the \emph{critical radius}.
Typically an upper bound on the localized
Gaussian width would allow this critical radius to serve
as an upper bound for the error $\ltwo{\genparamhat - \genparamstar}$.
Here, we show that with two-sided control
of the localized Gaussian width
and a regularity assumption,
the error also satisfies a matching lower bound.
In the next section, we will illustrate the consequence of this
result with some examples.


\subsection{Adaptive estimation rates}
\label{SecExamples}

We now demonstrate the consequences of \autoref{PropEllipseEstimation}
via some examples.  We begin with the simple problem of estimation for
$\genparamstar = 0$, where we see a number of standard rates from the
ellipse estimation literature.  We then consider some more interesting
examples of extremal vectors, and show how the resulting estimation rates
differ from the classical ones.


\subsubsection{Estimating at \texorpdfstring{$\genparamstar = 0$}{}}

We begin our exploration by considering ellipse-constrained
estimation problem at $\thetastar = 0$.
In this section, we focus on two type of ellipses that are specified
by aspect ratios $\eig{j}$ where $\eig{j}$ follows an
$\smoothpar$-polynomial decay and $\gamma$-exponential
decay. The first one corresponds to estimating a
function in $\smoothpar$-smooth Sobolev class---that is, functions
that are almost everywhere $\smoothpar$-times differentiable,
and with the derivative $f^{(\smoothpar)}$ being Lebesgue integrable.


\paragraph{$\alpha$-polynomial decay:}

Consider an ellipse $\Ellipse$ defined by the aspect ratios $\eig{j} =
c j^{-2 \smoothpar}$ for some $\smoothpar > 1/2$.  In
\autoref{ExampleSobolev}, inequality~\eqref{EqnRegAtZero}, it is
verified that this ellipse is regular at $0$, and that
\mbox{$\kind(\delta) \asymp \delta^{-1/\smoothpar}$}.  Thus, solving
the fixed point problem~\eqref{EqnFixedPt} yields $\delcrit \asymp
\noisestd^{\frac{2 \smoothpar}{2 \smoothpar + 1}}$, and one can check
that the conditions for \autoref{PropEllipseEstimation} are met.  Here
our notation $\asymp$ denotes equality up to constants independent of
$(\sigma, \usedim)$.  With a rescaling argument~\eqref{EqnRescaling},
the proposition implies
\begin{align}
c \left(\noisestd^2 \right)^{\frac{2
\smoothpar}{2 \smoothpar+1}} \,\le\, \ltwo{\genparamhat -
\genparamstar}^2 \,\le\, C \left(\noisestd^2 \right)^{\frac{2 \smoothpar}{2
\smoothpar+1}},
\end{align}
with probability $\ge 1 - \exp\parens*{- c' \noisestd^{-\frac{2}{2
\smoothpar + 1}}}$ for some constants $C > c > 0$ and $c'$.
One may notice that the rate
$\left(\noisestd^2 \right)^{\frac{2 \smoothpar}{2\smoothpar+1}}$ coincides
with the minimax estimation rate for estimating in an $\alpha$-smooth Sobolev
function class. We will show in our later section that it is indeed the case.


\paragraph{$\gamma$-exponential decay:}
Consider another case where the ellipse $\Ellipse$ is defined
by the aspect ratios $\eig{j} = c_1 \exp(-c_2 j^\gamma)$, for some $\gamma > 1/2$.
Then a slight modification of the computation in \autoref{ExampleSobolev}
yields
\begin{align*}
   \kind(\delta) = \argmin_k \{\sqrt{\eig{k+1}}, \frac{9}{10}\delta\} \asymp \log^{\frac{1}{\gamma}} \parens*{\frac{1}{\delta}}.
\end{align*}
In order to establish the regularity condition,
notice that in this case, $\inf_{\gamma \in \Gammaclass(\genparamstar, \delta)} \sum_{i=1}^\usedim \gamma_i$ is achieved in limit by 
$\gamma_i = \mu_{i} \indic\{i > \lwk(\delta)\}$ and further more
\begin{equation}
\label{Eqnchipmunk}
\inf_{\gamma \in \Gammaclass(\genparamstar, \delta)}
\sum_{i=1}^\usedim \gamma_i
= \sum_{j = \kind+1}^\usedim c_1 e^{-c_2 j^\gamma}
\asymp \int_{\kind}^\infty e^{-c_2 t^{\gamma}} \mathop{dt}
\le \frac{1}{\gamma \kind^{\gamma - 1}} \int_{\kind^\gamma}^\infty e^{-c_2 u} \mathop{du}
\asymp \frac{\eig{\kind}}{\kind^{\gamma - 1}}
\asymp \frac{\delta^2}{\kind^{\gamma - 1}},
\end{equation}
which by definition, shows that $\Ellipse$ is regular at $\genparamstar=0$.

Solving the fixed point problem~\eqref{EqnFixedPt} yields
$\delcrit \asymp \noisestd \log^{\frac{1}{2\gamma}} \parens*{\frac{1}{\noisestd}}$
up to other polylogarithmic factors in $\noisestd$.
One can check that the conditions for \autoref{PropEllipseEstimation}
are met, so by the rescaling argument~\eqref{EqnRescaling}, we have,
up to polylogarithmic factors,
\begin{align}
c \noisestd^2 \log^{\frac{1}{\gamma}} (\noisestd^{-1})
\,\le\,
\ltwo{\genparamhat - \genparamstar}^2
\,\le\, C \noisestd^2 \log^{\frac{1}{\gamma}} (\noisestd^{-1}),
\end{align}
with probability
$\ge 1 - \exp{- c' \log^{\frac{1}{\gamma}}(1 / \noisestd)}$
for some constants $C > c > 0$ and $c'$.



\subsubsection{Estimating at extremal vectors}
\label{SecExtremal}

In the previous section, we studied the adaptive estimation rate for $\at = 0$.
In this section, we study some non-zero cases of the vector $\at$.
For concreteness, we restrict our
attention to vectors that are non-zero some coordinate $s \in
[\usedim] = \{1, \ldots, \usedim\}$, and zero in all other
coordinates.  Even for such simple vectors, our analysis reveals some
interesting and adaptive scalings.

\begin{figure}[H]
\begin{center}
  \widgraph{0.7\textwidth}{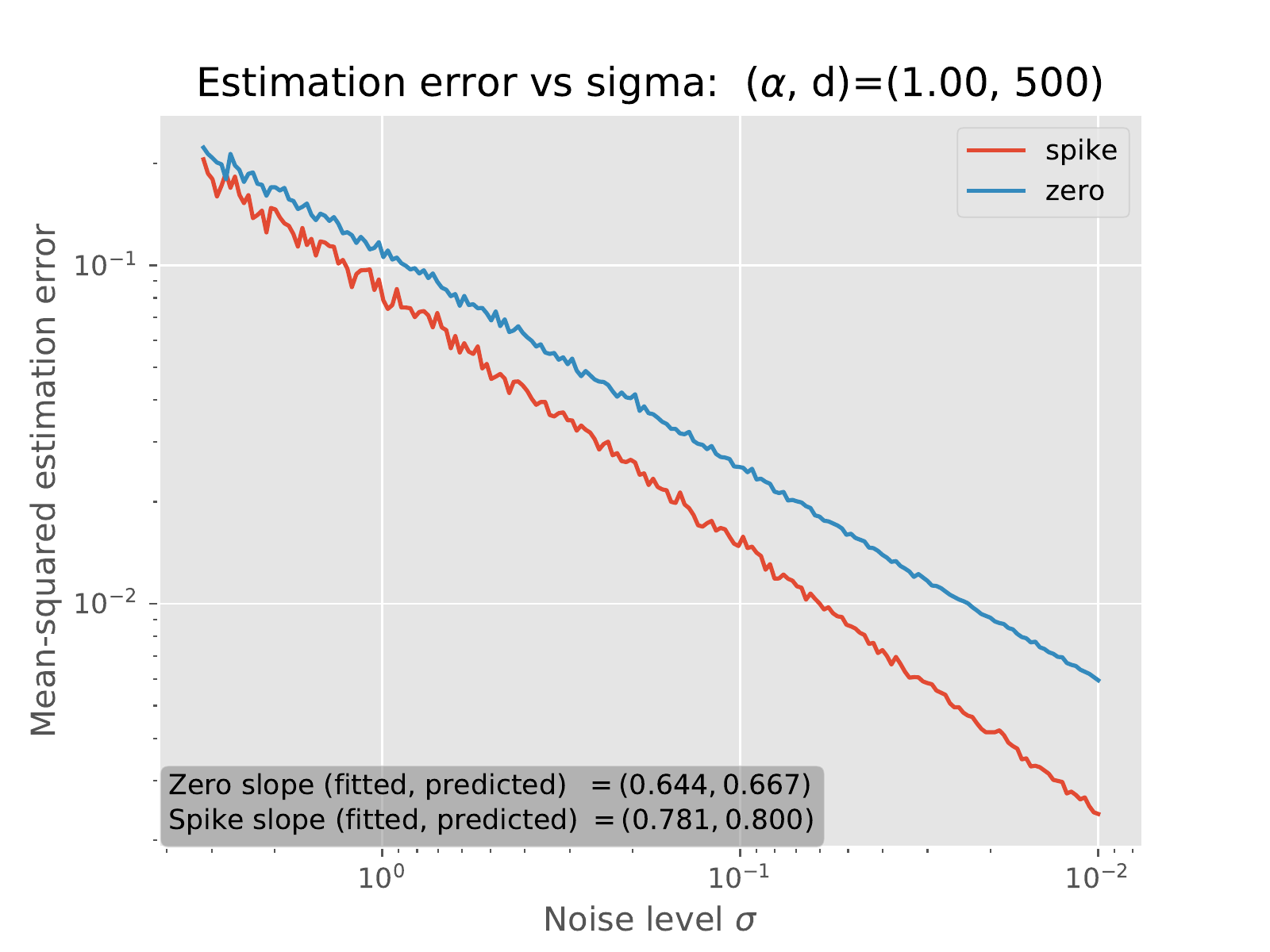}
\caption{Log-log plot of $\Exs \ltwo{\genparamhat - \genparamstar}^2$
  against $\noisestd$ for the ellipse with polynomial decay $\eig{j} =
  j^{-2}$ in $\usedim = 500$ dimensions. The blue curve is the case
  $\genparamstar = 0$, and the red curve is $\genparamstar = e_1$. }
\label{FigPoly}
\end{center}
\end{figure}

Given integer $s \in [\usedim]$, consider $\genparamstar \defn
(\sqrt{\eig{s}} - r) e_s$ for some $r \in [\newepscritl(s, \ellip),
  ~\newepscritu(s, \ellip)]$ where $\newepscritl(s, \ellip),
\newepscritu(s, \ellip)$ are small constants that are defined in
\citet[Corollary 2]{wei2017testing}.  Note that the shrinkage $-r$
away from the boundary is due to the boundary issue in
\autoref{ThmLowerBound}.  We believe it is an artifact of our analysis
that is possibly removable; for instance, in our simulations below
(\autoref{FigPoly}) we have an example with $\genparamstar =
\sqrt{\mu_1} e_1$ on the boundary of the ellipse that exhibits the
same predicted behavior as its shrunken counterpart.

So as to streamline notation, we adopt $\kind(\delta)$ as a short hand 
for $\kind(\at, \delta)$.
\citet{wei2017testing} (Section 4.4) show that with
$\xi = (1-\tinyconst)\delta$, we have
\begin{align}
\label{EqnSalt}
   \kind(\delta) = \kind \Big(\frac{\xi}{1-\tinyconst} \Big) \leq
   \underbrace{\arg \max_{1 \leq k \leq \usedim} \Big \{\mu_k^2 \geq
     \frac{1}{64}\xi^2\mu_s \Big \}}_{ = : \; m_u}.
\end{align}
This upper bound is proved by considering the projection onto
the $m_u$-dimensional subspace spanned by $\{e_1,\ldots,e_{m_u}\}$.
At the same time, we prove in \autoref{LemKolLinf} that
\begin{align}
\label{EqnPepper}
 \kind(\delta) \geq 0.09 \cdot m_{\ell}, \qquad \mbox{where $m_{\ell}
   \defn \arg \max \limits_{1 \leq k \leq \usedim} \Big \{\mu^2_k \geq
   \delta^2 \mu_s \Big \}$.}
\end{align}


\paragraph{$\alpha$-polynomial decay:}
Consider an ellipse $\Ellipse$ with $\eig{j} = c j^{-2 \smoothpar}$
for some $\smoothpar > 1/2$.  From the above calculation, we can
conclude that
\begin{align*}
    m_u,~m_\ell,~\upk \asymp (\mu_s \delta^2)^{-\frac{1}{4\alpha}},
\end{align*}
Here our notation $\asymp$ denotes equality up to constants
independent of problem parameters such as $(\sigma, \usedim)$.
Let us verify the regularity condition~\eqref{EqnRegularity} with dimension $m_u$ and projection
to linear space $\Proj_{m_u}$ spanned by $\{e_1,\ldots,e_{m_u}\}$.
Since $\gamma = \delta^2(0,\ldots,0,{\mu_{m_u+1}},\ldots {\mu_d})$ is feasible
in limit for the set $\Gammaclass(\genparamstar, \delta, m_u, \Proj_{m_u})$, we have
\begin{align*}
\inf_{\gamma \in \Gammaclass(\genparamstar, \delta, m_u, \Proj_{m_u})} \gamma_i
\leq
\delta^2 \sum_{i=m_u+1}^\usedim {\mu_{i}} \, \asymp \, \delta^2 \int_{m_u+1} t^{-2\alpha} dt
= \delta^2 (m_u+1)^{-2\alpha+1}.
\end{align*}
Since $\alpha > 0$, $\upk$ and $m_u$ is equal up to a constant,
the right hand side above is bounded above by $\delta^2\upk$, which establishes
the regularity condition at $\at.$

As long as $s \lesssim (\sigma^2)^{-2/(4\alpha+1)}$, solving the fixed point problem~\eqref{EqnFixedPt} yields $\delcrit \asymp
\noisestd^{\frac{4 \smoothpar}{4 \smoothpar + 1}}$,
and one can check that the conditions for \autoref{PropEllipseEstimation}
are met. Thus,
\begin{align}
c \left(\noisestd^2 \right)^{\frac{4\smoothpar}{4\smoothpar+1}} 
\, \le\, \ltwo{\genparamhat -
\genparamstar}^2 \,\le\, C \left(\noisestd^2 \right)^{\frac{4\smoothpar}{4\smoothpar+1}},
\end{align}
with probability $\ge 1 - \exp\parens*{- c' \noisestd^{-\frac{2}{4
      \smoothpar + 1}}}$ for some constants $C > c > 0$ and $c'$.


\paragraph{$\gamma$-exponential decay:}
Now consider ellipse $\Ellipse$ with $\eig{j} = c_1 \exp(-c_2 j^\gamma)$ for some $\gamma > 1/2$.
From the above calculation, we can conclude that
\begin{align*}
    m_u,~m_\ell,~\upk \asymp \log^{\frac{1}{\gamma}} \parens*{\frac{1}{\delta}}.
\end{align*}
Let us verify the regularity condition~\eqref{EqnRegularity} with dimension $m_u$ and
projection to linear space $\Proj_{m_u}$ spanned by $\{e_1,\ldots,e_{m_u}\}$.
Since $\gamma = \delta^2(0,\ldots,0,\sqrt{\mu_{m_u+1}},\ldots \sqrt{\mu_d})$ is feasible
for the set $\Gammaclass(\genparamstar, \delta, m_u, \Proj_{m_u})$, by similar calculation
from inequality \eqref{Eqnchipmunk}, we can show that the ellipse is regular at $\at$.

Solving the fixed point problem~\eqref{EqnFixedPt} yields
$\delcrit \asymp \noisestd \log^{\frac{1}{2\gamma}} \parens*{\frac{1}{\noisestd}}$
up to other polylogarithmic factors in $\noisestd$.
One can check that the conditions for \autoref{PropEllipseEstimation}
are met, so by the rescaling argument~\eqref{EqnRescaling}, we have,
up to polylogarithmic factors,
\begin{align}
c \noisestd^2 \log^{\frac{1}{\gamma}} (\noisestd^{-1})
\,\le\,
\ltwo{\genparamhat - \genparamstar}^2
\,\le\, C \noisestd^2 \log^{\frac{1}{\gamma}} (\noisestd^{-1}),
\end{align}
with probability
$\ge 1 - \exp\parens*{- c' \log^{\frac{1}{\gamma}}(\noisestd^{-1})}$
for some constants $C > c > 0$ and $c'$.


\paragraph{Numerical results:}
To illustrate our findings from above, \autoref{FigPoly} provides a
numerical plot of the mean-squared error of the constrained least
squared estimator~\eqref{EqnLS} for estimating the vector $\thetastar
= 0$ (blue curve) and the vector $\thetastar = e_1$ (red curve).  In
each case, the plot shows show the error decreases as a function of
the inverse noise level $\frac{1}{\sigma^2}$.

The underlying ellipse is defined by the eigenvalues $\eig{j} = j^{-2
  \smoothpar}$ with $\smoothpar = 1$.  Consequently, the predicted
scaling of the mean=squared error is
$(\noisestd^2)^{\frac{2\smoothpar}{2\smoothpar + 1}}$ for the zero
vector, and $(\noisestd^2)^{\frac{4\smoothpar}{4\smoothpar + 1}}$ for
the ``spiked'' $e_1$ vector.  Based on these predictions, our our
theory suggests that on a log-log plot, the mean-squared error should
decay at a linear rate with slopes $-2/3$ and $-4/5$ respectively.
The empirical least-squares fit shows that these predictions are very
accurate.


\subsection{Minimax risk bounds}
\label{SecMinimax}

As another consequence of our main results, in this section, we show
that the LSE is minimax optimal for ellipse estimation problem that is
described above.  Here the \emph{minimax risk} over the ellipse
$\ellip$ is defined as
\begin{align}
\MiniMax(\Ellipse) \defn \inf_{\genparamhat} \sup_{\genparamstar
  \in \ellip} \Exs_\genparamstar \ltwo{\genparamhat -
  \genparamstar}^2,
\end{align}
where the supremum is taken over distributions
$\NORMAL(\genparamstar, \noisestd^2 \IdMat_n)$
indexed by $\genparamstar \in \Ellipse$,
and the infimum is taken over all estimators.
By this criteria, estimators are compared on their worst-case
performance.

In the following, we show that the minimax optimal risk is achieved
by the LSE estimator and the risk is characterized through
the solution to the fixed point problem~\eqref{EqnFixedPt}.
Let $\delcrit(0)$ be the solution to the
fixed point problem~\eqref{EqnFixedPt} for $\at = 0$.

\begin{cors}
\label{CorMinimax}
There are universal constants $c,C > 0$, 
the global minimax risk of estimation over the entire ellipse
$\Ellipse$ satisfies
\begin{subequations}
\begin{align}
\label{EqnMinimaxLower}
\MiniMax(\Ellipse) \,\ge\, c\noisestd^2 \lwk(0, \delcrit).
\end{align}
If furthermore the ellipse is regular~\eqref{EqnRegularity} for all
$\genparamstar \in \Ellipse$, then
\begin{align}
\label{EqnMinimaxUpper}
\MiniMax(\Ellipse) \,\le\, C\noisestd^2
\lwk(0, \frac{1}{2} \delcrit).
\end{align}
\end{subequations}
\end{cors}
\noindent We prove this result in \autoref{SecMinimaxProof}.

In contrast to the minimax lower bound of \citet{yang2015randomized},
our minimax lower bound~\eqref{EqnMinimaxLower} does not require the
regularity assumption~\eqref{EqnRegularity}.  
See \autoref{SecKernelRegularity} for a discussion of how the notion of
regularity of \citet{yang2015randomized} is a special case of our
notion.  The lower bound is proved by showing that the ellipse
contains a $\kind$-dimensional ball, and then applying the standard
minimax bound in for estimation in a $\kind$-dimensional space.

On the other hand, the upper bound~\eqref{EqnMinimaxUpper} does
require the regularity assumption, which allows us to apply
\autoref{PropEllipseEstimation}. It implies that the risk of the LSE
for each problem $\genparamstar \in \Ellipse$ is upper bounded by
$\lesssim \delcrit^2(\genparamstar)$.  Furthermore, we show that among
all $\genparamstar$, the largest upper bound
$\delcrit^2(\genparamstar)$ is the case $\genparamstar = 0$, which
yields the upper bound in \autoref{CorMinimax}.  Thus, the hardest
problem for the LSE is estimating $\at = 0$, and its risk there
matches the lower bound.  In short, the LSE is minimax optimal for
ellipses that are regular.

\section{Proofs}
\label{SecProofs}

We now turn to the proofs of our main results, namely
\autoref{ThmUpperBound} and \autoref{ThmLowerBound}.  The proofs of
more technical results are deferred to appendices, as noted within
this section.


\subsection{Proof of \autoref{ThmUpperBound}}
\label{SecPfUB}

For any dimension and projection pair $(k, ~\Proj_k)$, we can write
\begin{align}
\label{EqnTwoTerms}
\Exs \sup_{\Delta \in \Ellipse_{\at} \cap \Ball(\delta)}
\inprod{\sgauss}{\Delta} & \leq \underbrace{\Exs \sup_{\Delta \in
    \Ellipse_{\at} \cap \Ball(\delta)}
  \inprod{\sgauss}{\Proj_{k}\Delta}}_{\Term_1} + \underbrace{\Exs
  \sup_{\Delta \in \Ellipse_{\at} \cap \Ball(\delta)}
  \inprod{\sgauss}{\Delta - \Proj_{k} \Delta}}_{\Term_2}.
\end{align}
We now proceed to upper bound the two terms $\Term_1$ and $\Term_2$.

\paragraph{Bounding $\Term_1$:}
From standard properties of orthogonal projections onto subspaces, we
have $\inprod{\sgauss - \Proj_{k} \sgauss}{\Proj_{k} \Delta} =
0$ for any $\sgauss$ and $\Delta$.  By combining this fact with the
Cauchy-Schwarz inequality, the term $\Term_1$ is upper bounded as
\begin{subequations}
\begin{align}
\Term_1 \; = \; \Exs \sup_{\Delta \in \Ellipse_{\at} \cap
  \Ball(\delta)} \inprod{\sgauss}{\Proj_{k}\Delta} = \Exs
\sup_{\Delta \in \Ellipse_{\at} \cap \Ball(\delta)}
\inprod{\Proj_{k} \sgauss}{\Proj_{k}\Delta} \le \Exs
\sup_{\Delta \in \Ellipse_{\at} \cap \Ball(\delta)}
\vecnorm{\Proj_{k} \sgauss}{2} \vecnorm{\Proj_{k}\Delta}{2}.
\end{align}
By the non-expansiveness of projection onto a subspace, we have
$\ltwo{\Proj_{k} \Delta} \le \ltwo{\Delta} \overset{(i)}{\leq}
\delta$, where inequality (i) follows from the inclusion $\Delta \in
\Ball(\delta)$.  Thus, we have established that
\begin{align}
\label{EqnTerm1Bound}
\Term_1 \leq \delta \Exs \vecnorm{\Proj_{k} \sgauss}{2} \leq \delta
\sqrt{k},
\end{align}
where the last step follows from first applying Jensen's inequality,
and then noting that the distribution of $\Proj_{k} w$ is a
$k$-dimensional standard Gaussian vector.

\paragraph{Bounding $\Term_2$:}
For a given vector $\gamma \in \Gammaclass(\genparamstar, \delta)$,
define the diagonal matrix $\Amat \defn \diag(\sqrt{\gamma_1}, \ldots,
\sqrt{\gamma_d})$.  Noting that $\inprod{\sgauss}{\Delta -
  \Proj_{k} \Delta} = \inprod{\Amat \sgauss}{\Amat^{-1} (\Delta -
  \Proj_{k} \Delta)}$ and then applying the Cauchy-Schwarz
inequality, we find that
\begin{align}
\Term_2 & \leq \Exs \sup_{\Delta \in \Ellipse_{\at} \cap
  \Ball(\delta)} \vecnorm{\Amat \sgauss}{2} \; \vecnorm{\Amat^{-1}
  (\Delta - \Proj_{k} \Delta)}{2}.
\end{align}
By the definition of $\Gammaclass(\genparamstar, \delta, k, \Proj_k)$,
we must have $\vecnorm{\Amat^{-1} (\Delta - \Proj_{k} \Delta)}{2} \le
1$.  Thus, we have the upper bound
\begin{align}
\label{EqnTerm2Bound}
\Term_2 & \leq \Exs \ltwo{\Amat \sgauss} \le \sqrt{\sum_{i=1}^\usedim
  \gamma_i},
\end{align}
where the last step is due to Jensen's inequality.
\end{subequations}
Since our choice of $\gamma$ was arbitrary, we may add an infimum over
$\gamma \in \Gammaclass(\genparamstar, \delta, k, \Proj_k)$.
Combining the two bounds~\eqref{EqnTerm1Bound}
and~\eqref{EqnTerm2Bound} concludes the proof.

\subsection{Proof of \autoref{ThmLowerBound} }
\label{SecPfLB}

As in the preceding proof, we adopt $\lwk$ as convenient shorthand for
the quantity $\lwk(\genparamstar, \delta)$.
We now divide our analysis into two cases, depending on
whether or not $\enorm{\at} \leq 1/2$.


\subsubsection{Case I}

First, suppose that $\enorm{\at} \leq \frac{1}{2}$, which implies that
$\Rfun(\delta) \leq (\enorm{\at}-1)^2\leq 1$.  Under this condition,
Lemma 2 from the paper~\cite{wei2017testing} guarantees that
\begin{align}
\label{EqnByEndOfToday}
\kwidth_{k}(\ellip_{\thetastar}\cap \Ball((1 - \tinyconst) \delta)) \,
\leq\, \frac{3}{2} \min \Big \{ (1 - \tinyconst) \delta,
\sqrt{\eig{k+1}} \Big\}.
\end{align}
By definition, the critical dimension $\lwk \defn \arg \min \limits_{k=1,
  \ldots, \usedim} \{\kwidth_{k}(\Ellipse_{\at} \cap \Ball((1 - \tinyconst) \delta))
\leq \frac{9}{10} \delta\}$ can be upper bounded as
\begin{align}
  \label{EqnDefnLwkPrime}
\lwk \, \leq \, \arg \min_{k = 1, \ldots, \usedim} \{ \frac{3}{2}
\sqrt{\eig{k+1}} \leq \frac{9}{10} \delta\} \, = \,: \, \lwk',
\end{align}
where we have used the fact that $\frac{9}{10} \leq 1 - \tinyconst$, and
$\kwidth_{k}(\Ellipse_{\at} \cap \Ball((1 - \tinyconst) \delta))$ is non-decreasing
in $k$.

Let $E_{\lwk'}$ denotes the $\lwk'$-dimensional subspace of  vectors that are zero is their last $\usedim - \lwk'$ coordinates.
Recalling that $\Sph(r)$ denotes a Euclidean sphere of radius $r$, we
claim that
\begin{align}
\label{EqnChain_Case1}
  \GWidth(\Ellipse_{\at} \cap \Ball(\delta)) \, \overset{(i)}{\geq}
  \, \GWidth \Big( \Sph(\frac{3}{10}\delta) \cap E_{\lwk'}
  \Big) \, \overset{(ii)}{=} \, \frac{3}{10} \delta \sqrt{\lwk'}.
\end{align}
Taking this claim as given for the moment, combining it with the
bounds $\enorm{\at} \leq 1/2$ and $\lwk \leq \lwk'$, we find that
\begin{align*}
  \GWidth(\Ellipse_{\at} \cap \Ball(\delta)) \geq \frac{3}{10} \delta
  \sqrt{1 - \enorm{\genparamstar}^2} \sqrt{\lwk(\genparamstar,
    \delta)},
\end{align*}
which completes the proof of \autoref{ThmLowerBound} in this case.


\paragraph*{Proof of inequality~\eqref{EqnChain_Case1}:}

In this proof, we adopt the convenient shorthand $b = 3/10$.  Part
(ii) of the inequality can be seen from the spherical example in the
discussion of \autoref{ThmUpperBound}.  It only remains to prove part
(i).  Let us first show that $\Sph(2 b\delta) \cap E_{\usedim - \lwk'}
\subset \ellip$. Recalling the definition of $\lwk'$ from
equation~\eqref{EqnDefnLwkPrime}, we have
\begin{align*}
\sum_{i=1}^{\usedim} \frac{x_i^2}{\eig{i}} = \sum_{i=1}^{\lwk'}
\frac{x_i^2}{\eig{i}} \, \overset{(iii)}{\leq} \, \sum_{i=1}^{\lwk'}
\frac{x_i^2}{\eig{\lwk'}} = \frac{(2b\delta)^2}{\eig{\lwk'}} \,
\overset{(iv)}{\leq} \, 1,
\end{align*}
where inequality (iii) follows from the non-increasing order of
$\eig{i}$ and inequality (iv) follows from the definition of $\lwk'$.

In order to establish the inclusion $\Ball_{\lwk'}(b\delta) \subset
\ellip_{\at}$, we make use of the fact that $\enorm{\at} \leq 1/2$.
Since $\enorm{2\at} \leq 1$, we have $2\at \in \ellip$. For any $v \in
\Sph_{\lwk'}(b\delta)$, since $\Ball_{\lwk'}(2b\delta) \subset \ellip$
we have $2v \in \ellip$.  Combining these two facts together and the
convexity of set $\ellip$, we have $v+\at \in \ellip$.  It further
implies that $\Ball_{\lwk'}(b\delta) \subset \ellip_{\at}$ and
finishes the proof of inequality~\eqref{EqnChain_Case1}.


\subsubsection{Case II}

Otherwise, we may assume that $\enorm{\at} > 1/2$, in which case
$\Rfun(\delta/\ccon) \leq (\enorm{\at}-1)^2 < 1$, and hence by
definition of the function $\Rfun$, we have $\delta <
\ccon\ltwo{\at}/a$.  For the remainder of the proof, we assume that
$\lwk \geq 160$.  The case when $\lwk < 160$ is addressed separately
at the end of this proof.

The proof of \autoref{ThmLowerBound} requires two auxiliary lemmas.
The first is packing lemma, proved in \citet[Lem. 4]{wei2017testing}.
Here we state a slightly altered version of this claim, better suited
to our purposes.  Let $\NewMat$ denote the diagonal matrix with
entries $1/\mu_1, \ldots, 1/\mu_d$, and adopt the shorthands $a \defn
1 - \tinyconst$ and $b \defn \frac{3}{10}$ based on the definition of
the critical dimension~\eqref{Eqnkcrit}.
\begin{lems}
\label{LemPacking}
For any vector $\thetastar \in \ellip$ such that $\ltwo{\thetastar} >
a \epsilon$, there exists a vector $\thetadag \in \Ellipse$, a
collection of $\usedim$-dimensional orthonormal vectors $\{u_i
\}_{i=1}^{\lwk}$ and an upper triangular matrix of the form
\begin{align}
\Mat \defn
    \begin{bmatrix}
        1 & h_{3,2} & h_{4,2} & \cdots & h_{\lwk,2}\\ & 1 & h_{4,3} &
        \cdots & h_{\lwk,3}\\ &&1 & \cdots & h_{\lwk,4} \\ &&&\ddots &
        \vdots \\ &&&&1
    \end{bmatrix} \in \real^{\lwk - 1, \lwk - 1}
\end{align}
with ordered singular values $\nu_1 \ge \cdots \ge \nu_{\lwk-1} \ge 0$
such that:
\begin{enumerate}[(a)]
\item The vectors $u_1$, $\NewMat \thetadag$, and $\thetadag -
  \thetastar$ are all scalar multiples of one another.
\item We have $\ltwo{\thetadag - \thetastar} = a \delta$.
\item
  Letting $H_{\cdot, i}$ denote the $i$th column of $\Mat$, for every
  $i \in [\lwk - 1]$, the vector $\thetadag \pm b \delta
  \underbrace{\begin{bmatrix} u_2 & \cdots u_{\lwk}
\end{bmatrix}}_{\defn U} H_{\cdot, i}$
belongs to the ellipse $\Ellipse$.

  \item We have $\enorm{\thetadag} \le \enorm{\genparamstar}$.

\item For any integers $t_1 \in [\lwk-1]$, $t_2 \in [\lwk-2]$, we have
  \begin{align}
    \label{EqnEignBound}
\nu_{t_1} \overset{(i)}{\leq} \frac{a}{3b}\sqrt{\frac{\lwk-1}{t_1}},
~~\text{ and }~~ \nu_{t_2+1} \overset{(ii)}{\geq} 1 -
\frac{t_2}{\lwk-1} - \sqrt{\frac{a^2 - 9b^2}{9b^2}}.
\end{align}
\end{enumerate}
\end{lems}

Before proving \autoref{ThmLowerBound}, let us introduce some
notation.  Let $\Mat$, $U$ and $\thetadag$ be as given in the
\autoref{LemPacking} above and let $X \defn U \Mat$ have columns
$x_1,\ldots, x_{\kind - 1}$ Let $V$ be the matrix of right singular
vectors of $\Mat$ so that $\Mat^\top \Mat = V \Sigma^2 V^\top$, where
$\Sigma^2$ is diagonal with the squared singular values $\nu_1^2 \ge
\cdots \ge \nu_{\lwk - 1}^2$ of $\Mat$ in order.




Let $m_1 \defn \floor{(\lwk - 1) / 8}$ and $m_2 \defn \floor{(\lwk -
  1) / 4}$, and define the sparsity level $s \defn \rho \frac{\kind -
  1}{16}$ for some constant\footnotemark{} $\rho \in (0,1)$.
\footnotetext{The arguments that follow do not depend on the specific
  choice of $\rho$, and taking $\rho = 1/2$ suffices.  However in the
  proof of \autoref{CorMetEnt}, we re-use these arguments for a
  different value of $\rho$.}  For a given $s$-sized subset $\Sset$ of
$\{m_1, \ldots, m_2 \}$, any vector of the form $\zs = (\zs_1, \ldots,
\zs_{\lwk - 1}) \in \{-1, 0, 1\}^{\lwk - 1}$ with zeros in all
positions not indexed by $\Sset$ is called as an \emph{$\Sset$-valid
  sign vector}.  Any such sign vector can be used to define the
perturbed vector
 \begin{align}
\label{EqnPackElt}
 \genparam^S \defn \thetadag + b \delta \frac{1}{\sqrt{32 s}} U \Mat V z^\Sset
 \end{align}
The following lemma guarantees the existence of a large collection
$\Tclass$ of $s$-sized subsets of $\{m_1, \ldots, m_2 \}$ such that
the collection $\{\genparam^S, S \in \Tclass \}$ has certain desirable
properties.
\begin{lems}
\label{LemPert}
There exists a collection $\Tclass$ of $s$-sized subsets of $\{m_1,
\ldots, m_2\}$ such that:
\begin{enumerate}[(a)]
  \item The collection $\Tclass$ has cardinality at least $\binom{
    \lfloor \frac{1}{16} (\lwk - 1) \rfloor }{s}$.
  \item For each $\Sset \in \Tclass$, there is a $\Sset$-valid sign
    vector $\zs$ such that the associated perturbation $\genparam^S$
    belongs to the ellipse $\Ellipse$, and moreover satisfies the
    bounds:
\begin{align}
  \label{EqnCond2}
  \delta^2 \; \overset{(i)}{\le} \; \vecnorm{\genparam^S -
    \genparamstar}{2}^2 \; \overset{(ii)}{\le} \;
  \frac{4}{1-\enorm{\genparamstar}^2} \: \delta^2.
\end{align}
\end{enumerate}
\end{lems}
\noindent See Appendix~\ref{SecProofSparseCombination} for the proof
of this lemma. \\

Turning back to the proof of \autoref{ThmLowerBound}, consider those
perturbation vectors \eqref{EqnPackElt} that are defined via
\autoref{LemPert}. For each $\Sset \in \Tclass$, we define the vectors
\begin{align*}
\Deltatilde^S \defn \genparam^S - \genparamstar, \quad \mbox{and}
\quad \Delta^S \defn \frac{\delta}{\vecnorm{\Deltatilde^S}{2}}
\Deltatilde^S.
\end{align*}
Inequality~\eqref{EqnCond2} implies that
$\frac{\delta}{\vecnorm{\Deltatilde^S}{2}} \leq 1$.  By the convexity
of the set $\Ellipse_\at$, we have $\Delta^S \in \Ellipse_{\at} \cap
\Sph(\delta)$ for each $S \in \Tclass$.  By restricting the supremum
to a smaller subset, we obtain the lower bound
\begin{align}
\Exs \sup_{\Delta \in \Ellipse_{\at} \cap \Sph(\delta)}
\inprod{\sgauss}{\Delta} \ge \Exs \max_{S \in \Tclass}
\inprod{\sgauss}{\Delta^S}.
\end{align}
Re-writing the definition~\eqref{EqnPackElt} in the form $\genparam^S
= \thetadag + \frac{b \delta}{\sqrt{32 s}} U \Mat V z^S$, it follows
that
\begin{align}
\Delta^S \defn \frac{\delta}{\ltwo{\Deltatilde^S}} \Deltatilde^S =
\delta \parens*{ \frac{1}{\ltwo{\Deltatilde^S}} (\thetadag -
  \genparamstar) + \frac{b \delta}{\sqrt{32 s} \ltwo{\Deltatilde^S}} U
  \Mat V z^S },
\end{align}
which further guarantees that
\begin{align}
\Exs \max_{\Sset \in \Tclass} \inprod{\sgauss}{\Delta^S} & \ge \delta
\Exs \max_{\Sset \in \Tclass}
\inprod{\sgauss}{\frac{b\delta}{\sqrt{32 s}\vecnorm{\Deltatilde^S}{2}} U
  \Mat V z^S} + \delta \Exs \max_{\Sset \in \Tclass}
\inprod{\sgauss}{\frac{1}{\vecnorm{\Deltatilde^S}{2}} (\thetadag -
  \genparam^*)} \\
& =\Exs \max_{\Sset \in \Tclass}
\inprod{\sgauss}{\frac{b\delta}{\sqrt{32 s}\vecnorm{\Deltatilde^S}{2}} U
  \Mat V z^S},
\end{align}
where the second equality follows since $\Exs
\inprod{\sgauss}{\thetadag - \genparamstar} = 0$.  The right-hand side
is non-negative, since for any fixed choice of $S_0 \in \Tclass$, we
have
\begin{align}
\Exs \max_{S \in \Tclass}
\inprod{\sgauss}{\frac{1}{\vecnorm{\Deltatilde^S}{2}} U \Mat V z^S}
\ge \Exs \inprod{\sgauss}{\frac{1}{\vecnorm{\Deltatilde^{S_0}}{2}} U
  \Mat V z^{S_0}} = 0.
\end{align}
Noting that inequality~\eqref{EqnCond2}(ii) can be rewritten as
$\ltwo{\Deltatilde^S}^2 \le \frac{4}{1 - \enorm{\genparamstar}^2}
\delta^2$, we find that
\begin{align}
\Exs \max_{S \in \Tclass}
\inprod{\sgauss}{\frac{1}{\vecnorm{\Deltatilde^S}{2}} U \Mat V z^S}
\ge \sqrt{\frac{1 - \enorm{\genparamstar}^2}{4 \delta^2}} \; \Exs
\max_{S \in \Tclass} \inprod{\sgauss}{U \Mat V z^S}.
\end{align}
Putting together the pieces, we have established that
\begin{align}
\label{EqnBigStep1}
\Exs \sup_{\Delta \in \Ellipse_{\at} \cap \Sph(\delta)}
\inprod{\sgauss}{\Delta} \ge \frac{b \delta}{16} \sqrt{\frac{1 -
    \enorm{\genparamstar}^2}{s}} \; \Exs \max_{S \in \Tclass}
\inprod{\sgauss}{U \Mat V z^S}.
\end{align}
Our next step is to lower bound the expected maximum on the RHS, and
to this end, we state an auxiliary result:
\begin{lems} Under the conditions  of \autoref{ThmLowerBound},
we have
\label{LemSudakov}
\begin{align}
\label{EqnSudakov}
  \Exs \max_{S \in \Tclass} \inprod{\sgauss}{U \Mat V z^S} \geq
   \frac{1}{4} \Exs \Big[\max_{S \in \Tclass} \sum_{i \in S} \sgauss_i \Big].
\end{align}
\end{lems}
\noindent See Appendix~\ref{SecPfLemSudakov} for the proof
of this lemma.

Let us now control the term on the right-hand side of
inequality~\eqref{EqnSudakov}.  Let $\Aevent$ be the event that there
are least $s$ positive elements among the i.i.d. standard Gaussian
random variables $\{ \sgauss_i \}_{i={m_1}}^{m_2}$.  By the law of
total expectation, we have
\begin{align}
\Exs \Big[ \max_{S \in \Tclass} \sum_{i \in S} \sgauss_i \Big] & =
\underbrace{\Exs \Big[ \max_{S \in \Tclass} \sum_{i \in S} \sgauss_i
    \; \mid \; \Aevent \Big]}_{\Term_1} \; \Prob[\Aevent] +
\underbrace{\Exs \Big[ \max_{S \in \Tclass} \sum_{i \in S} \sgauss_i
    \; \mid \; \Aevent^c \Big]}_{\Term_2} \; \Prob[\Aevent^c].
\end{align}
Beginning our analysis with $\Term_1$, under the event $\Aevent$,
there exists some (random) subset $S' \in \Tclass$ of cardinality
$|S'| \geq s$ such that $\sgauss_i > 0$ for all $i \in S'$.  (When
there are multiple such sets, we choose one of them uniformly at
random.)  In terms of this set, we have
\begin{align}
\Term_1 \; = \; \Exs \Big[ \max_{S \in \Tclass} \sum_{i \in S}
  \sgauss_i \; \mid \; \Aevent \Big] & \ge \Exs_{w, S'} \Big[ \sum_{i
    \in S'} \sgauss_i \; \mid \; \Aevent \Big] \; = \; \sum_{S'}
\Exs_{w} \Big[ \sum_{i \in S'} \sgauss_i \mid S' \Big] \; \mprob[S'
  \mid \Aevent],
\end{align}
where $\mprob[S' \mid \Aevent]$ denotes the conditional probability of
the randomly chosen $S'$ given that $\Aevent$ holds.  Since we are
conditioning on a random set $S'$ on which each $\sgauss_i$ is
positive, we have
\begin{align*}
\Exs_{w} \Big[ \sum_{i \in S'} \sgauss_i \mid S' \Big] & = \Exs_{w}
\Big[ \sum_{i \in S'} \sgauss_i \mid \sgauss_i > 0 \Big] \\
& \geq s \, \Exs[\sgauss_i \mid \sgauss_i > 0] \; = \; s \sqrt{2/\pi}.
\end{align*}
Since $\sum_{S'} \mprob[S' \mid \Aevent] = 1$, we have proved that
$T_1 \geq s \sqrt{2 / \pi}$.

Turning to the term $\Term_2$, we begin by observing that for any
fixed $S_0 \in \Tclass$, we have
\begin{align}
\max_{S \in \Tclass} \sum_{i \in S} \sgauss_i \ge \sum_{i \in S_0}
\sgauss_i \ge - \sum_{i \in S_0} |\sgauss_i|.
\end{align}
Using this observation we can conclude that
\begin{align}
\Exs \Big[ \max_{S \in \Tclass} \sum_{i \in S} \sgauss_i \; \mid \;
  \Aevent^c \Big] \ge  \Exs\Big[ -\sum_{i \in S_0}|\sgauss_i|  \; \mid \;
  \Aevent^c \Big] \overset{(i)}{=}  \Exs \Big[- \sum_{i \in S_0} |\sgauss_i| \Big] = - s \sqrt{2 / \pi}.
\end{align}
where (i) follows from the fact $\Aevent^c$ only depends on the sign
of $w_i$ and the distribution of $|\sgauss_i|$ is independent of
$\Aevent^c$.  Combining these two lower bounds, we find that
\begin{align}
\Exs \max_{S \in \Tclass} \sum_{i \in S} \sgauss_i \ge s \sqrt{2 /
  \pi} (1 - 2 \Prob[\Aevent^c]).
\end{align}

We now bound the probability of event $\Aevent^c$.  Recall that event
$\Aevent$ holds if and only if there are at least $s$ positive
elements among the i.i.d. standard Gaussian random variables $\{
\sgauss_i \}_{i={m_1}}^{m_2}$.  Since $s \defn \floor{(m_2-m_1) / 4}$,
with probability no larger than $\exp(-
(m_2-m_1)\kullinf{\frac{1}{4}}{\frac{1}{2}}) \le e^{-0.1 (m_2-m_1)}$,
there are more than $(m_2-m_1)/ 4$ components among $\sgauss_{m_1},
\ldots, \sgauss_{m_2}$ that are positive, meaning that
$\Prob[\Aevent^c] \le e^{-0.1 (m_2-m_1)}$.  Thus, we have the lower
bound
\begin{align}
\Exs \Big[ \max_{S \in \Tclass} \sum_{i \in S} \sgauss_i \ge s \sqrt{2
    / \pi} (1 - 2 e^{-0.1 (m_2-m_1)}) \Big] & \geq \frac{1}{5} s,
\end{align}
where the last step uses the fact that $m_2-m_1 \ge \lwk/16 > 10$.

Combining this last bound with inequalities~\eqref{EqnBigStep1}
and~\eqref{EqnSudakov} yields
\begin{align}
\Exs \sup_{\Delta \in \Ellipse_{\at} \cap \Sph(\delta)}
\inprod{\sgauss}{\Delta} & \ge \frac{b \delta}{16} \sqrt{\frac{1 -
    \enorm{\genparamstar}^2}{s}} \; \; \Exs \max_{S \in \Tclass}
\inprod{\sgauss}{U \Mat V z^S} \\
& \geq \frac{b}{64} \delta \sqrt{\frac{1 -
    \enorm{\genparamstar}^2}{s}} \; \; \Exs \Big[\max_{S \in \Tclass}
  \sum_{i \in S} \sgauss_i \Big] \\
& \ge \frac{b}{320} \delta \sqrt{(1 - \enorm{\genparamstar}^2) s}
\\
& \ge c' \sqrt{1 - \enorm{\genparamstar}^2} \cdot \delta \sqrt{\lwk},
\end{align}
where the last step uses the fact that $s = \rho \frac{\kind - 1}{16}$.

In order to finish the proof, we deal with the case of $\lwk < 160$
separately.  According to part (b) of \autoref{LemPacking}, if we
denote $v_1 \defn \thetastar - \thetadag$, then $\thetadag \in
\Ellipse$ and $\ltwo{v_1}= a\delta.$ It is also shown in the proof of
\citet[Lem. 5]{wei2017testing} that $\thetastar + v_1 \in \Ellipse.$
Therefore the two points $\pm v_1$ are both contained in
$\Ellipse_{\at} \cap \Ball(\delta)$ for a sufficiently small $\delta$.
As a result, we have $\GWidth(\Ellipse_{\at} \cap \Sph(\delta)) \geq
\GWidth(\{\pm v_1 \}) = a \delta \sqrt{2/\pi}$, which establishes the
lower bound in \autoref{ThmLowerBound} with constant $c' = \frac{a}{4
  \sqrt{5\pi}}$.

\section{Discussion}

In this paper, we studied the behavior of localized Gaussian widths
over ellipses.  These localized widths are known to play a fundamental
role in controlling the difficulty of associated testing and
estimation problems.  Despite its fundamental importance, the
localized Gaussian width is hard to compute in general.  The main
contribution of our paper was to show how the localized Gaussian width
can be bounded, both from above and below, via the localized
Kolmogorov dimension.  These Kolmogorov dimensions can be computed in
many interesting cases, which leads to an explicit characterization of
the estimation error of least-squares regression as a function of the
true regression vector within the ellipse.  We used this
characterization to show how the difficulty of estimating a vector
$\thetastar$ within the ellipse can vary dramatically as a function of
the location of $\thetastar$.  Estimating the all-zeros vector
($\thetastar = 0$) is always the hardest sub-problem, and leads to the
global minimax rate.  Much faster rates of estimation can be obtained for
vectors located near ``narrower'' portions of the ellipse boundary.
While much of the analysis in this paper is specific to ellipses, we
do anticipate that the general procedure of moving from Gaussian width
to the Kolmogorov width could be useful in studying adaptivity and
local geometry in other estimation problems.

\subsection*{Acknowledgments}
This work was partially supported by Office of Naval Research grant
DOD-ONR-N00014, and National Science Foundation grant NSF-DMS-1612948
to MJW.  In addition, BF was partially supported by a National Science
Foundation Graduate Research Fellowship.

\vspace*{2cm}


\appendix


\section{Properties of kernel regularity}
\label{SecKernelRegularity}

In this section, we relate our definition of
regularity~\eqref{EqnRegularity} to a concept introduced in previous
work by \citet{yang2015randomized}.  In the context of kernel ridge
regression, they defined the quantity
\begin{align}
\label{EqnYangStatDim}
\upktilde \equiv \upktilde(\delta) \defn \arg\min_k \, \{\eig{k+1} \le
\delta^2\}
\end{align}
with the convention $\upktilde = \usedim$ if the minimization is over
an empty set. They said that an ellipse is \emph{regular} if
\begin{align}
\label{EqnYangRegularity}
\sum_{j = \upktilde + 1}^\usedim \eig{j} \le c \upktilde \delta^2,
\qquad \text{for all $\delta > 0$,}
\end{align}
where $c>0$ is some universal constant that does not depend on
$\usedim$.  They used this property to prove a minimax lower bound on
the prediction error for kernel ridge regression.

Let us now show that our regularity assumption~\eqref{EqnRegularity}
is a generalization of the condition~\eqref{EqnYangRegularity}, in
that it reduces to it in the special case $\genparamstar = 0$.  In
order to establish this claim, we begin by observing that for any $k
\in \{1, \ldots, \usedim - 1\}$, we have $\kwidth_{k}(\Ellipse_{\at}
\cap \Ball((1 - \tinyconst) \delta)) = \min\{\eig{k+1}^{1/2}, (1 - \tinyconst) \delta\}$ because the
minimization in the definition~\eqref{EqnKolmogorov} is achieved by
the projection onto the subspace $\myspan\{e_1, \ldots, e_k\}$, and
the maximization is achieved by $\theta = \min\{\eig{k+1}^{1/2},
(1 - \tinyconst) \delta\} e_{k+1}$.  On the other hand, for $k = \usedim$ we have
$\kwidth_{k}(\Ellipse_{\at} \cap \Ball((1 - \tinyconst) \delta)) = 0$.  Putting these
two together gives
\begin{align}
\upk(0,\delta) = \min \Big \{ k \mid \eig{k+1} \le \frac{81}{100}
\delta^2 \Big\},
\end{align}
(with the convention $\upk = \usedim$ if the minimum is over an empty
set). Thus, we have recovered definition~\eqref{EqnYangStatDim} up to
a constant factor in $\delta$.

Since the optimal projection $\Proj_{\upk}$ is the projection onto the
linear subspace $\myspan\{e_1, \ldots, e_{\upk}\}$, we can consider a
sequence of positive vectors approaching $\gamma \defn (\eig{i}
\Indic{i > k})_{i=1}^d$ to obtain
\begin{align*}
\inf_{\gamma \in \Gammaclass(\genparamstar,
  \delta)}\sum_{i=1}^{\usedim}\gamma_i \le \sum_{i = \upk + 1}^\usedim
\eig{i}.
\end{align*}
Consequently, our regularity condition~\eqref{EqnRegularity} holds as
long as $\sum_{i = \upk + 1}^\usedim \eig{i} \le c \upk \delta^2$.
Thus, it matches the notion of regularity~\eqref{EqnYangRegularity}
considered in~\citet{yang2015randomized}.


\section{Proof of \autoref{CorMetEnt}}
\label{AppProofCorMetEnt}

Throught out this proof, we use $c,c',c''$ etc. to denote universal
constants that do not depend on any problem parameters such as
$\delta, \mu_i$ and $\thetastar$ and their values can vary from line
to line.

The proof of inequality (i) in equation~\eqref{EqnMetEnt} is
straightforward.  By combining the Sudakov
minoration~\eqref{EqnSudakovMinoration} with our upper
bound~\eqref{EqnGWidthSandwich} on the localized Gaussian width, we
find that
\begin{align*}
c' \delta \sqrt{\log \Mpack(\delta / 2, \Ellipse_{\at} \cap
  \Ball(\delta))} \le \GWidth(\Ellipse_{\at} \cap \Ball(\delta)) \le
\constup \delta \sqrt{\upk(\genparamstar, \delta)}.
\end{align*}
Thus, we have proved inequality (i) in equation~\eqref{EqnMetEnt}.

We now turn to the proof the second inequality (ii).  It is convenient
to divide our analysis into two cases depending on whether or not
$\enorm{\genparamstar} \le \frac{1}{2}$.


\paragraph*{Case 1:} $\enorm{\genparamstar} \le \frac{1}{2}$.
As shown earlier in equation~\eqref{EqnChain_Case1} from the proof of
\autoref{ThmLowerBound}, the set \mbox{$\Ellipse_{\at} \cap
  \Ball(\delta)$} contains the $\kind'$-dimensional sphere
$\Sph(\frac{3}{10} \delta) \cap E_{\kind'}$.  Thus, by a standard
volume argument~\cite{pisier1986probabilistic,Wai17}, it must have log
packing number bounded from below by $c \kind' \log \frac{1}{\delta}$.
This quantity is lower bounded by $\kind$ up to some universal
constant, which establishes inequality (ii) in this case.


\paragraph*{Case 2:} $\enorm{\genparamstar} > \frac{1}{2}$.
We follow the notation from \autoref{SecProofs}.  In the proof of
\autoref{ThmLowerBound} (in particular, see
equation~\eqref{EqnPackElt} and \autoref{LemPert}), we constructed a
set of vectors $\genparam^\Sset$ that after rescaling, all lie in our
set $\Ellipse_{\at} \cap \Ball(\delta)$.  Each such vector
$\genparam^\Sset$ is formed by taking a certain point $\thetadag$ near
$\genparamstar$, and adding certain combinations of orthogonal vectors
$u_i$.  We argue here that there is a subset of these scaled vectors
of size $\gtrsim \kind$ that are pairwise separated from each other by
a distance $\gtrsim \delta$.

We are only interested in proving bounds up to constant factors,
meaning that we may assume without loss of generality that $\kind \ge
32 \times 10^4$; otherwise the result~\eqref{EqnMetEnt} holds
immediately with a sufficiently large choice of $c'$.

Recall the earlier definition $s \defn \rho \frac{\kind - 1}{16}$ for
a fixed constant $\rho \in (0, 1)$; for this argument, we take $\rho =
10^{-4}$.  By Lemma 4.10 in \citet{Massart03}, we can find a
subset of $s$-sparse vectors contained in the binary hypercube $\{0,
1\}^{\frac{1}{16}(\kind - 1)}$ with log cardinality at least
\begin{align}
s \log \frac{\frac{1}{16}(\upk-1)}{s} \gtrsim \lwk,
\end{align}
and such that any pair of distinct elements differs in at least $(2 -
2\rho)s$ entries.  Transferring this result to the context of
\autoref{LemPert}, we are guaranteed a collection of vectors of log
cardinality $\gtrsim \lwk$ such that
\begin{align}
\label{EqnSignPacking}
\ltwo{z^\Sset - z^{\Sset'}}^2 > (2 - 2 \rho) s
\end{align}
for $z^\Sset \ne z^{\Sset'}$ in our packing.

Recalling that $V^\top \Mat^\top \Mat V = \Sigma^2$ and the
definition~\eqref{EqnPackElt} of $\genparam^\Sset$, we then have
\begin{align}
\ltwo{\genparam^\Sset - \genparam^{\Sset'}}^2 &= \frac{b^2
  \delta^2}{32 s} \ltwo{U \Mat V (z^\Sset - z^{\Sset'})}^2 \\ &=
\frac{b^2 \delta^2}{32 s} (z^\Sset - z^{\Sset'})^\top \Sigma^2
(z^\Sset - z^{\Sset'}).
\end{align}
Since $z^\Sset$ and $z^{\Sset'}$ are zero in their first $m_1 - 1$
components, we can use inequality (ii) from \autoref{LemPacking} to
bound the relevant diagonal entries of $\Sigma$. Doing so yields
\begin{align}
\ltwo{\genparam^\Sset - \genparam^{\Sset'}}^2 & \ge \frac{b^2
  \delta^2}{32 s} \parens*{1 - \frac{1}{8} - \sqrt{\frac{a^2 - 9
      b^2}{9 b^2}}}^2 \ltwo{z^\Sset - z^{\Sset'}}^2.
\end{align}
Thus, we have obtained a collection of vectors $\genparam^\Sset$,
indexed by subsets $\Sset$, such that \mbox{$\ltwo{\genparam^\Sset -
    \genparam^{\Sset'}} \gtrsim \delta$} for $\Sset \ne \Sset'$.

Finally, we need to show that after shrinking these $\genparam^\Sset$
toward $\genparamstar$ and re-centering, we obtain a packing of
$\Ellipse_{\at} \cap \Ball(\delta)$.  For each $\Sset$ recall the
definitions $\Deltatilde^\Sset \defn \genparam^\Sset - \genparamstar$
and $\Delta^\Sset \defn \frac{\delta}{\vecnorm{\Deltatilde^\Sset}{2}}
\Deltatilde^\Sset$.  From discussion below \autoref{LemPert}, we have
already showed that each vector $\Delta^\Sset$ lies in $\Ellipse_{\at}
\cap \Ball(\delta)$; it only remains to verify that distinct pairs
are well-separated.

First, direct computation yields
\begin{align}
  \label{EqnPackInprod}
\ltwo{\Delta^\Sset - \Delta^{\Sset'}}^2 &= 2 \delta^2 \Big (1 -
\frac{\inprod{\Deltatilde^{\Sset}}{\Deltatilde^{\Sset'}}}
     {\ltwo{\Deltatilde^{\Sset}}\ltwo{\Deltatilde^{\Sset'}}} \Big).
\end{align}
In order to show that the right-hand side is lower bounded by a
constant multiple of $\delta^2$, it suffices to upper bound the inner
product term.  Using the fact that $\thetadag - \genparamstar$ has
norm $a \delta$ and is orthogonal to the columns of $U$ (see
\autoref{LemPacking}), we have
\begin{align}
 \inprod{\Deltatilde^S}{\Deltatilde^{S'}} &= \inprod{\thetadag -
   \genparamstar + \frac{b\delta}{\sqrt{32 s}} U\Mat V z^\Sset}
        {\thetadag - \genparamstar + \frac{b\delta}{\sqrt{32 s}} U\Mat
          V z^{\Sset'}} \\
& = a^2 \delta^2 + \frac{b^2 \delta^2}{32 s} z^\Sset \Sigma^2
        z^{\Sset'}.
\end{align}
If $z^\Sset \ne z^{\Sset'}$ are from our packing, then by construction
they differ on at least $(2 - 2 \rho) s$ components, so they must
agree on at most $\rho s$ components. Applying the inequality (i) from
\autoref{LemPacking} to bound the relevant entries of $\Sigma^2$, we
can continue from above to obtain
\begin{align}
\inprod{\Deltatilde^S}{\Deltatilde^{S'}}
    &\le a^2 \delta^2 + \frac{b^2 \delta^2}{32} \cdot
    8 (\frac{a}{3b})^2 \cdot \rho
    \\
    &\le a^2 \parens*{1 + \frac{\rho}{36}} \delta^2 < \delta^2.
\end{align}
The last inequality follows from our earlier choice of $a \defn 1 -
10^{-5}$ and $\rho \defn 10^{-4}$.  Dividing both sides by
$\ltwo{\Deltatilde^\Sset} \ltwo{\Deltatilde^{\Sset'}} \ge \delta^2$
(where this inequality follows from \autoref{LemPert}), we can
continue from our earlier step~\eqref{EqnPackInprod} to obtain
\begin{align}
\ltwo{\Delta^\Sset - \Delta^{\Sset'}}^2 \ge c \delta^2.
\end{align}
Putting together the pieces, we have exhibited the claimed packing of
$\Ellipse_{\at} \cap \Ball(\delta)$ of log cardinality $\gtrsim \kind$
and packing radius $\gtrsim \delta$.


\section{Proof of \autoref{CorMinimax}}
\label{SecMinimaxProof}

We divide our proof into two parts, corresponding to the upper and
lower bounds respectively.

\paragraph*{Upper bound:}
Let us start with the proof of the upper bound.  Under the regularity
assumption, we may apply \autoref{PropEllipseEstimation} to bound the
mean-squared error $\Exs_{\at} \ltwo{\genparamhat - \genparamstar}^2$
of the LSE; in particular, it is upper bounded by $\delcrit^2(\at)$ up
to an universal constant.  (Recall that $\delcrit(\at)$ is the
solution to the fixed point equation~\eqref{EqnFixedPt}.)

In order to arrive at the desired minimax upper bound, we need to show
that the function $\at \mapsto \delcrit(\at)$ is maximized at $\at =
0$.  Since $\kind$ is a non-increasing function of $\delta$ (see the
paper~\cite[Sec. D.1]{wei2017testing}), a larger $\kind(\at)$
corresponds to a larger value of $\delcrit(\at)$.  These two
quantities are related via the equation
\begin{align*}
\delcrit(\at) &= \constlw \noisestd \sqrt{\lwk(\at, \delcrit)}.
\end{align*}

\noindent The following lemma bounds the supremum of $\lwk$.
\begin{lems}
\label{Lemfun}
The critical dimensions at any $\at$ can be controlled as 
\begin{align}
\upk(\genparamstar, \delta) \leq \kind(0, \frac{1}{2} \delta) +
1 \qquad \mbox{for all $\delta \in \Big(0, \;
  \Phi^{-1}((\enorm{\at}^{-1}-1)^2) \wedge \sqrt{\mu_1} \Big)$.}
\end{align}
\end{lems}
\noindent The proof of this lemma is given in
Appendix~\ref{SecProofLemfun}.  Note that it implies the claimed upper
bound upper bound~\eqref{EqnMinimaxUpper}.


\paragraph*{Lower bound:}

By definition, the minimax risk decreases when the supremum is taken
over a smaller subset. In order to establish the lower bound, we
restrict the supremum to a ball around zero.  Recall our calculations
from \autoref{ExampleSobolev}, where we showed that the Kolmogorov
width of a local ball around $\thetastar = 0$ is given by
\begin{align*}
  \kwidth_{k}(\Ellipse_{\at} \cap \Ball((1 - \tinyconst) \delta) =
  \min \big \{ \sqrt{\mu_{k+1}}, (1 - \tinyconst) \delta \big \}.
\end{align*}
The corresponding $\lwk(0, \delta)$ is given by
\begin{align*}
\lwk(0, \delta) \defn \arg \min_{k = 1, \ldots, \usedim} \Big \{
\kwidth_{k}\Big(\Ellipse \cap \Ball\big((1 -
\tinyconst)\delta\big)\Big) \leq \frac{9}{10} \delta \Big \}.
\end{align*}
By inspection, we have the upper bound $\lwk(0, \delta) = \arg \min_{k
  = 1, \ldots, \usedim} \big \{ \sqrt{\mu_{k+1}} \leq \frac{9}{10}
\delta \Big \}$.  We also have the lower bound $\sqrt{\mu_{\lwk(0,
    \delta)}} \geq \frac{9}{10} \delta$ for every $\delta \leq
\sqrt{\mu_1}$.  Note that the ellipse $\Ellipse$ always contains a
$k$-dimensional ball centered at zero with radius $\sqrt{\mu_k}$.
Combined with the bounds just stated, for every $\delta \in (0,
\sqrt{\mu_1}]$, the ellipse also contains a ball of radius
  $\frac{9}{10} \delta$ centered at zero of dimension $\lwk(0,
  \delta)$.

Now we are ready to control the minimax risk.  First notice that
\begin{align}
\label{EqnBanff}
\MiniMax(\Ellipse) \defn \inf_{\genparamhat} \sup_{\at \in \ellip}
\Exs_\at \ltwo{\genparamhat - \at}^2 \geq \inf_{\genparamhat}
\sup_{\at \in \Ball(\frac{9}{10} \delta) \cap E_{\kind(0,\delta)}}
\Exs_\genparam \ltwo{\genparamhat - \at}^2,
\end{align}
where recall that $E_{m}$ denotes the space which contains
$d$-dimensional vectors with their last $d-m$ coordinates all equal to
zero.

By standard results (e.g., see the book~\cite{Wai17}), estimating a
$m$-dimensional vector in a $r$ radius ball has minimax risk lower
bounded as
\begin{align*}
  \inf_{\genparamhat} \sup_{\at \in \Ball(r) \cap E_{m}}
  \Exs_\genparam \ltwo{\genparamhat - \at}^2
  \gtrsim \min \{r^2, m \noisestd^2 \}.
\end{align*}
Substituting this lower bound into inequality~\eqref{EqnBanff}, we
find that
\begin{align}
\label{EqnBanff2}
\MiniMax(\Ellipse) \gtrsim \min \{ \big(\frac{9}{10}\delta\big)^2,
\kind(0,\delta) \noisestd^2 \},
\end{align}
for each $\delta \leq \sqrt{\mu_1}$.  From the
definition~\eqref{EqnFixedPt}, we have $\delcrit(0) = \constlw
\noisestd \sqrt{\lwk(0, \delcrit)}$.  Taking $\delta = \delcrit(0)$ in
inequality~\eqref{EqnBanff2} yields the claimed lower
bound~\eqref{EqnMinimaxLower}.




\section{Proof of \autoref{PropEllipseEstimation}}
\label{SecEstimationProof}

This appendix is devoted to the proof of
\autoref{PropEllipseEstimation}.

\subsection{Reduction to bounding localized Gaussian width}
\label{SecReduction}

\citet{chatterjee2014new} provided one way of obtaining upper and
lower bounds on the error $\ltwo{\genparamhat - \genparamstar}$ of
the least squares estimator for a general convex set, under the Gaussian
sequence model~\eqref{EqnGaussianModel}.  Define the function
\begin{align}
\label{EqnCritFun}
\CritFun(t) \defn \frac{\delta^2}{2} - \noisestd
\GWidth(\Ellipse_{\at} \cap \Ball(\delta)),
\end{align}
which can be shown to be strongly convex on $(0,\infty)$ with a unique
minimizer $\deltachat > 0$.  Then:
\begin{theos}[{\cite[Thm. 1.1, Cor. 1.2]{chatterjee2014new}}]
\label{ThmChatterjee}
The least squares estimator $\genparamhat$ satisfies
\begin{align}
\label{EqnChatterjeeOrig}
\big|\ltwo{\genparamhat - \genparamstar} - \deltachat\big| \le t
\sqrt{\deltachat}, \quad \text{w.p. } \ge 1 -
3\exp\parens*{-\frac{t^4}{32 \noisestd^2 (1 +
    t/\sqrt{\deltachat})^2}},
\end{align}
for any $t > 0$.  Furthermore, there is a universal constant $C > 0$
such that
\begin{align}
\label{EqnChatExpectationBound}
\big|\Exs \ltwo{\genparamhat - \genparamstar}^2 - \deltachat^2 \big|
\le C \deltachat^{3/2} \noisestd^{1/2}, \qquad \text{if $\deltachat
  \ge \noisestd$}.
\end{align}
\end{theos}

In particular, if we take $t = c\sqrt{\deltachat}$, it is guaranteed
that
\begin{equation}
\label{EqnChatterjee}
\big|\ltwo{\genparamhat - \genparamstar} - \deltachat\big|
\le c\deltachat,
\quad \text{w.p. } \ge 1 - 3\exp(-c' \deltachat^2 / \noisestd^2).
\end{equation}
The following simple lemma shows how sandwiching $\CritFun$ between
two functions allows us to obtain upper and lower bounds for its
minimizer $\deltachat$.

\begin{lems}
  \label{LemSandwich}
Suppose that there are functions $\CritFunLw, \CritFunUp$ such that
$\CritFunLw(\delta) \le \CritFun(\delta) \le \CritFunUp(\delta)$ for
all $\delta \in [0, \infty)$.  Then for any $r \geq \inf
  \limits_{\delta \geq 0} \CritFunUp(\delta)$, we have
\begin{align}
\deltachat \in \{\delta \ge 0 : \CritFunLw(\delta) \le r\}.
\end{align}
In particular, if $\CritFunLw$ is unimodal, then this sub-level set is
an interval.
\end{lems}
The proof of this lemma is simple.  For a given $r \geq
\inf_{\delta \geq 0} \CritFunUp(\delta)$, we have
\begin{align*}
\CritFunLw(\deltachat) \overset{(i)}{\le} \CritFun(\deltachat)
\overset{(ii)}{=} \inf_{\delta \geq 0} \CritFun(\delta)
\overset{(iii)}{\le} \inf_{\delta \geq 0} \CritFunUp(\delta) \; \leq
\; r
\end{align*}
where inequalities (i) and (iii) follow from the assumed sandwich
relation, and equality (ii) follows from the fact that $\deltachat$ is
the minimizer of $\CritFun$.

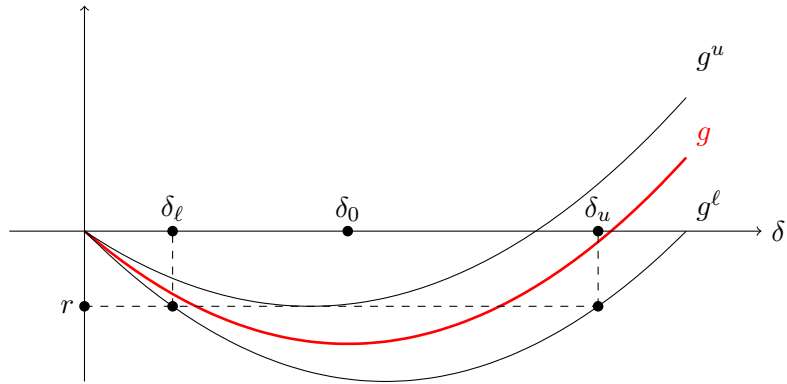
\begin{figure}[H]
\centering
\begin{tikzpicture}[scale=2]
\draw[domain=0:4, smooth, variable=\x, red, line width = 1pt] plot ({\x}, {12/49 * (\x-7/4)^2 - 3/4}); 
\draw[domain=0:4, smooth, variable=\x] plot ({\x}, {2/9 * (\x-3/2)^2 - 1/2}); 
\draw[domain=0:4, smooth, variable=\x] plot ({\x}, {1/4 * (\x-2)^2 - 1}); 
\draw[->] (-0.5,0) -- (4.5,0) node[right] {$\delta$};
\draw[->] (0,-1) -- (0,1.5) node[above] {};
\node[left] at (0, -1/2) {$r$};
\fill (0, -1/2) circle[radius=1pt];
\draw[dashed] (0, -1/2) -- (3.4142, -1/2);
\draw[dashed] (0.5858, -1/2) -- (0.5858, 0);
\draw[dashed] (3.4142, -1/2) -- (3.4142, 0);
‎
\node[above] at (7/4,0) {$\deltachat$};
\fill (7/4,0) circle[radius=1pt];
\fill (0.5858,-1/2) circle[radius=1pt];
\fill (3.4142, -1/2) circle[radius=1pt];
\node[above] at (0.5858,0) {$\delta_\ell$};
\fill (0.5858,0) circle[radius=1pt];
\node[above] at (3.4142,0) {$\delta_u$};
\fill (3.4142, 0) circle[radius=1pt];
\node[above right] at (4,0) {$\CritFunLw$};
\node[above right, red] at (4,1/2) {$\CritFun$};
\node[above right] at (4,1) {$\CritFunUp$};
\end{tikzpicture}
\caption{ Visualization of \autoref{LemSandwich} when $r =
  \inf_{\delta \geq 0} \CritFunUp(\delta)$, and $\CritFunLw$ is
  convex.  }
\end{figure}

\autoref{LemSandwich} and the bound~\eqref{EqnChatterjee} together
show that bounds on the localized Gaussian width that appears in the
definition~\eqref{EqnCritFun} of $\CritFun$ can be used to obtain high
probability upper and lower bounds on the error of the LSE.

We remark that the case for estimation over $\Ellipse(\bigrad)$ for
$\bigrad > 0$ reduces to the case $\bigrad = 1$ by rescaling.  Let
$\Ellipse(\bigrad)_{\genparamstar} \defn \{\genparam - \genparamstar :
\genparam \in \Ellipse(\bigrad)\}$ denote the re-centered ellipse.
Note that $\CritFun$ can be rewritten as
\begin{align}
  \label{EqnRescaling}
\CritFun(\delta) \defn \frac{\delta^2}{2} - \noisestd
\GWidth(\Ellipse(\bigrad)_{\genparamstar} \cap \Ball(\delta)) = R^2
\underbrace{ \left[\frac{\deltatilde^2}{2} - \widetilde{\noisestd}
    \GWidth(\Ellipse_{\widetilde{\genparam}^*} \cap
    \Ball(\deltatilde)) \right] }_{\CritFunTilde(\deltatilde)}
\end{align}
after the changes of variables $\deltatilde \defn \delta / \bigrad$,
$\widetilde{\genparam}^* \defn \genparamstar / \bigrad$, and
$\widetilde{\noisestd} \defn \noisestd / \bigrad$.  Then one can focus
on bounding $\CritFunTilde$ and ultimately re-scale by $R$ any bounds
obtained for the minimizer of $\CritFunTilde$ in order to obtain
bounds for the original minimizer $\deltachat$.


\subsection{Main portion of the proof}
\label{SecPfProp1}

Recall the two functions defined in equation~\eqref{EqnCritBounds}.
Under our assumptions, the bounds~\eqref{EqnGWidthSandwich} hold, so
that the critical function $\CritFun$ from equation~\eqref{EqnCritFun}
is sandwiched as $\CritFunLw(\delta) \le \CritFun(\delta) \le
\CritFunUp(\delta)$ for all $\delta$.  Now \autoref{LemSandwich} is
applicable for ellipse $\Ellipse$, constant $r = -
\frac{\delcrit^2}{2} = \CritFunUp(\delcrit)$ and function pair
$(\CritFunUp, \CritFunLw)$, so we know $\deltachat \in \{\delta \ge 0
\; \mid \; \CritFunLw(\delta) \le r\}$.

Since function $\CritFunLw$ is convex in $\delta$, there are two
solutions $\delcritleft$ and $\delcritright$ to the equation
\begin{align}
\label{EqnGetLevelSet}
\CritFunLw(\delta) = - \frac{\delcrit^2}{2},
\end{align}
and \autoref{LemSandwich} guarantees that
\begin{align}
\label{EqnChatInterval}
\delcritleft \leq \deltachat \leq \delcritright.
\end{align}

Moreover, we show below that $c_1 \delcrit \le \deltachat \le c_2
\delcrit$.  Taking this inequality to be true for the moment,
combining it with equation~\eqref{EqnChatterjee}
yields
\begin{align}
(1 - c) c_1 \delcrit \le \ltwo{\genparamhat - \genparamstar} &\le
  (1+c) c_2 \delcrit, \qquad \text{with prob. } \ge 1 - 3 \exp(-c'
  \delcrit^2 / \noisestd^2),
\end{align}
which concludes the proof.
Note that we arrive at the expectation bounds~\eqref{EqnExpectationBound}
by simply applying the earlier result~\eqref{EqnChatExpectationBound}.

It remains to show that $c_1 \delcrit \leq \delcritleft$ and
$\delcritright \le c_2 \delcrit$.  After some manipulation using the
fixed point equation~\eqref{EqnFixedPt},
equation~\eqref{EqnGetLevelSet} can be rewritten as
\begin{align}
\parens*{\delta - \noisestd \sqrt{\constup^2 \upk(\delta)}}^2 = \noisestd^2
(\constup^2 \upk(\delta) - \constlw^2 \lwk(\delcrit)).
\end{align}
Note that the solutions $\delta$ to the
equality~\eqref{EqnGetLevelSet} must satisfy $\constup^2 \upk(\delta)
\ge \constlw^2 \lwk(\delcrit)$, as required for the right-hand side to
be non-negative. In addition, they must satisfy one of the following
two equations:
\begin{subequations}
\begin{align}
\delta &= \noisestd \sqrt{\constup^2 \upk(\delta)} + \noisestd \sqrt{\constup^2
  \upk(\delta) - \constlw^2 \lwk(\delcrit)} \defn h_+(\delta),
\label{EqnSubFixedPt1}
\\
\delta &= \noisestd \sqrt{\constup^2 \upk(\delta)} - \noisestd
\sqrt{\constup^2 \upk(\delta) - \constlw^2 \lwk(\delcrit)} \defn h_-(\delta).
\label{EqnSubFixedPt2}
\end{align}
\end{subequations}
Note that any solution $\delcritright$ to the first
equation~\eqref{EqnSubFixedPt1} is larger than any solution
$\delcritleft$ to the second equation~\eqref{EqnSubFixedPt2}.  Indeed,
we have $\delcritleft = h_-(\delcritleft) < h_+(\delcritleft)$, so the
non-increasing nature of $h_+$ guarantees that the solution
$\delcritright$ to the equation $\delta = h_+(\delta)$ must be larger
than $\delcritleft$.

\begin{itemize}
  \item We first consider the solution $\delcritright$ to the first
equation~\eqref{EqnSubFixedPt1}. It is easy to check that
\begin{align*}
  \noisestd \constup \sqrt{\upk(\delta)} \leq h_+(\delta) \leq 2\noisestd \constup \sqrt{\upk(\delta)}.
\end{align*}
Recall $\upk(\delta)$ is non-increasing in $\delta$.
We know $\delcritright$ is smaller than the
solution to $\delta = 2 \noisestd \constup \sqrt{\upk(\delta)}$,
which in turn is smaller than $c_2 \delcrit$
(by assumption (c) of \autoref{PropEllipseEstimation}).
We thus have $\delcrit \leq \delcritright \le c_2 \delcrit$.

\item Next we consider the solution $\delcritleft$ to the second
equation~\eqref{EqnSubFixedPt2}.
We claim that $\delcritleft \geq c_1 \delcrit$.
In order to show this, we prove that $h_-(\delta)$ satisfies
\begin{align}
\label{EqnHfunc}
  h_-(c_1 \delcrit) \overset{(i)}{\geq} c_1 \delcrit,
  \text{ for some } c_1 \in (0,1)
  \qquad
  \text{and }~~
  h_-(\delcrit) \overset{(ii)}{\leq} \delcrit.
\end{align}
Take the above inequalities as given for now, we can combine them
with the fact that $h_-(\delta)$ is a non-decreasing function of $\delta$
to conclude that
the fixed point solution $\delcritleft$ of \eqref{EqnSubFixedPt1} satisfies
$c_1 \delcrit  \leq \delcritleft \leq \delcrit$.
\end{itemize}
Putting these two pieces together with inequality~\eqref{EqnChatInterval},
we conclude the proof of \autoref{PropEllipseEstimation}.
It remains to prove the inequalities~\eqref{EqnHfunc}.

\paragraph*{Proof of part (i):}
Applying the simple inequality $\constup^2 \upk(c_1\delcrit) - \constlw^2
\lwk(\delcrit) \leq \constup^2 \upk(c_1\delcrit)$ yields
\begin{align*}
  h_-(c_1 \delcrit) = \frac{\noisestd \constlw^2 \lwk(\delcrit)}
  {\sqrt{\constup^2 \upk(c_1\delcrit)} + \sqrt{\constup^2 \upk(c_1\delcrit) -
      \constlw^2 \lwk(\delcrit)}} \geq \frac{\noisestd \constlw^2
    \upk(\delcrit)}{ 2\sqrt{\constup^2 \upk(c_1 \delcrit)}} \geq c_1
  \noisestd \sqrt{\constlw^2 \upk(\delcrit)},
\end{align*}
where the last inequality follows by the fact that $\constup^2\upk(c_1
\delta) \leq \frac{1}{4 c^2_1} \constlw^2\lwk(\delcrit)$ (cf. Assumption
(b) in \autoref{PropEllipseEstimation}).  The fixed point
equation~\eqref{EqnFixedPt} further implies that
\begin{align*}
 h_-(c_1 \delcrit) \geq c_1 \noisestd \sqrt{\constlw^2 \upk(\delcrit)} =
 c_1 \delcrit,
\end{align*}
which proves our claim (i).

\paragraph*{Proof of part (ii):} Using the fact that
$\constup^2 \upk(\delcrit) \geq \constlw^2 \lwk(\delcrit)$, we find that
\begin{align}
h_-(\delcrit) = \frac{\noisestd \constlw^2 \lwk(\delcrit)} {\sqrt{\constup^2
    \upk(\delcrit)} + \sqrt{\constup^2 \upk(\delcrit) - \constlw^2
    \lwk(\delcrit)}} \leq \frac{\noisestd \constlw^2
  \lwk(\delcrit)}{\sqrt{\constup^2 \upk(\delcrit)}} \leq \noisestd
\sqrt{\constlw^2 \lwk(\delcrit)} = \delcrit,
\end{align}
where the last equality follows from the fact that $\delcrit$ is a
solution of the fixed point equation.  This completes the proof of
claim (ii).


\section{Auxiliary proofs for \autoref{ThmLowerBound}}

In this appendix, we collect the proofs of various auxiliary results
that underlie \autoref{ThmLowerBound}.


\subsection{Proof of \autoref{LemPert}}
\label{SecProofSparseCombination}

The set class $\Tclass$ to be demonstrated consists of all $s$-sized
subsets of a particular subset \mbox{$\Tset \subset \{m_1, \ldots,
  m_2\}$;} the subset $\Tset$ is constructed to have cardinality at
least $\lfloor \frac{\lwk - 1}{16} \rfloor$, so that the set class
$\Tclass$ has at least $\binom{\lfloor \frac{1}{16} (\lwk - 1)
  \rfloor}{s}$ elements.

Consider the $\lwk - 1$ diagonal elements of the matrix $V^\top X^\top
B X V$.  The sum of these diagonal elements is $\mytrace{V^\top X^\top
  B X V}$. Furthermore, the pigeonhole principle ensures that the
smallest $\frac{15}{16} (\lwk - 1)$ of the diagonal elements are each
at most
\begin{equation}
\label{EqnTraceBound}
\frac{16}{\lwk - 1} \mytrace{V^\top X^\top B X V} = \frac{16}{\lwk -
  1} \mytrace{X^\top B X} \le 16 \max_{i \le \lwk - 1} \enorm{x_i}^2.
\end{equation}
Let $\Tset$ be the indices of those $\frac{15}{16} (\lwk - 1)$
diagonal elements that are also in $\{m_1, \ldots, m_2\}$.  By
construction, we have $|\Tset| \ge m_2 - m_1 - \frac{1}{16} (\lwk - 1)
= \frac{\lwk - 1}{16}$, as desired.

Given the set class $\Tclass$ defined by the subset $\Tset$, we now
show that inequality (i) in equation~\eqref{EqnCond2} holds. Note that
\autoref{LemPacking} implies that any sign vector $z^S$ supported on
$S$ satisfies
\begin{align}
\vecnorm{\genparam^S - \genparamstar}{2}^2 = \vecnorm{\thetadag -
  \genparamstar}{2}^2 + \frac{b^2 \delta^2}{32s} \vecnorm{U \Mat V
  z^S}{2}^2.
\end{align}
Here the decomposition uses the fact that $\thetadag - \genparamstar$
is parallel to $u_1$, as guaranteed by part (a) of
\autoref{LemPacking}; this property ensures that $\thetadag -
\genparamstar$ is orthogonal to $u_2, \ldots, u_{\lwk}$.  Since the
columns of $U$ are orthogonal unit vectors, we have $\vecnorm{U \Mat V
  z^S}{2}^2 = \vecnorm{\Mat V z^S}{2}^2$.  Then recalling that $V^\top
\Mat^\top \Mat V = \Sigma^2$ is a diagonal matrix containing the
squared singular values of $\Mat$, we may use inequality (ii) in
\autoref{LemPacking} to obtain
\begin{align}
\label{EqnConstConstraint}
\vecnorm{\genparam^S - \genparamstar}{2}^2 &= a^2 \delta^2 + \frac{b^2
  \delta^2}{32 s} \vecnorm{\Mat V z^S}{2}^2 \\ &\ge \Big[a^2 +
  \frac{b^2}{32} \Big( 1 - \frac{m_2}{\lwk - 1} - \sqrt{\frac{a^2 - 9
      b^2}{9 b^2}} \Big)^2\Big] \delta^2 \\ &\ge \Big(a^2 + \frac{
  b^2}{2^9}\Big) \delta^2
\end{align}
where the last step follows from inequality \eqref{EqnSVbound}.  Here
let us take $\tinyconst$ small enough, for instance $10^{-5}$ such
that the right hand side above is greater than $\delta^2$.  (We have
made these choices of constants for the sake of convenience in the
proof, but note that other choices of these quantities are also
possible.)

Now, we prove that $\genparam^S \in \Ellipse$ and inequality (ii) in
equation~\eqref{EqnCond2} holds, in particular by using a
probabilistic argument.  Recall that $B :=
\operatorname{diag}(\eig{1}^{-1}, \ldots, \eig{\usedim}^{-1})$ so that
$\enorm{x}^2 = x^\top B x$.  For a given subset $S$, we specify a
random choice of $z^S$, in which for each $j \in S$, the value $z_j^S
\in \{-1, +1\}$ is an independent Rademacher variable.  Using this
random choice of $z^S$, we then let $\genparam^S$ be defined as in
equation~\eqref{EqnPackElt}, so that it is now a random vector.

Now part (a) of \autoref{LemPacking} guarantees that the vector $B
\thetadag$ is orthogonal to $u_2, \ldots, u_{\lwk}$.  As a
consequence, we have $\enorm{\genparam^S}^2 = (\thetadag)^\top B
\thetadag + \frac{b^2 \delta^2}{32 s} \enorm{X V z^S}^2$.

Let us focus on the expectation of the second term in the equation
above.  By the linearity and cyclic invariance properties of trace, we
have
\begin{align}
\Exs \enorm{X V z^S}^2 &= \Exs [(z^S)^\top V^\top X^\top B X V z^S]
\\
&= \mytrace{V^\top X^\top B X V \Exs[z^S (z^S)^\top]} \\ &=
\mytrace{V^\top X^\top B X V I_S},
\end{align}
where $I_S = \Exs[z^S (z^S)^\top]$ is the diagonal matrix whose $i$th
diagonal entry is $1$ if $i \in S$ and zero otherwise.  The last
expression is the sum of $s$ diagonal entries of $V^\top X^\top B X V$
indexed by elements of $\Tset$, so that our earlier
bound~\eqref{EqnTraceBound} implies that
\begin{equation}
\Exs \enorm{X V z^S}^2 \leq \; 16 \, s \; \max_{i \le \lwk - 1}
\enorm{x_i}^2.
\end{equation}
Letting $i^*$ denote the maximizer of the right-hand side, then
combining the previous few displays yields
\begin{align}
\Exs \enorm{\genparam^S}^2 &= (\thetadag)^\top B \thetadag + \frac{b^2
  \delta^2}{32 s} \Exs \enorm{X V z^S}^2 \; \leq \; (\thetadag)^\top B
\thetadag + \frac{b^2 \delta^2}{2} \enorm{x_{i^*}}^2.
\end{align}
Again using the fact that $B \thetadag$ and $x_{i^*}$ are orthogonal,
we have $\enorm{\thetadag + b \delta x_{i^*}}^2 = \enorm{\thetadag}^2
+ b^2 \delta^2 \enorm{x_{i^*}}^2$, and thus
\begin{align}
\label{EqnMarkovPart1}
\Exs \enorm{\genparam^S}^2 \le \frac{1}{2} \enorm{\thetadag}^2 +
\frac{1}{2} \enorm{\thetadag + b \delta x_{i^*}}^2 \le
\frac{\enorm{\genparamstar}^2 + 1}{2},
\end{align}
where the last step is due to the fact that $\enorm{\thetadag} \le
\enorm{\genparamstar}$ by construction, as well as $\enorm{\thetadag +
  b \delta x_{i^*}}^2 \le 1$, by claim (c) in \autoref{LemPacking}.

Similarly, part (a) of \autoref{LemPacking} implies the vector
$\thetadag - \genparamstar$ is orthogonal to all of the vectors $u_2,
\ldots, u_{\lwk}$, whence
\begin{align}
\vecnorm{\genparam^S - \genparamstar}{2}^2 = \vecnorm{\thetadag -
  \genparamstar}{2}^2 + \frac{b^2 \delta^2}{2 s} \ltwo{U \Mat V
  z^S}^2.
\end{align}
By properties of the trace along with the fact that $I_S = \Exs [z^S
  (z^S)^\top]$ is the diagonal matrix with $i$th diagonal entry equal
to $1$ if $i \in S$ and zero otherwise, we then have
\begin{align}
\Exs \vecnorm{\genparam^S - \genparamstar}{2}^2 &= \vecnorm{\thetadag
  - \genparamstar}{2}^2 + \frac{b^2 \delta^2}{32 s} V^\top \Mat^\top
\Mat V I_S.
\end{align}
By noting $V^\top \Mat^\top \Mat V = \Sigma^2$ is diagonal with the
squared singular values of $\Mat$ and applying the bound (i) from
\autoref{LemPacking}, we have
\begin{align}
\label{EqnMarkovPart2}
\Exs \vecnorm{\genparam^S - \genparamstar}{2}^2 &\le \Big(a^2 +
\frac{b^2}{32} \cdot \frac{a^2}{9b^2} \cdot \frac{\lwk - 1}{m_1} \Big)
\delta^2 \le \Big(1 + \frac{1}{36}\Big)a^2 \delta^2 < 2 \delta^2.
\end{align}

By a union bound and Markov's inequality, the two
inequalities~\eqref{EqnMarkovPart1} and~\eqref{EqnMarkovPart2} imply
\begin{align}
\mprob \Big( \enorm{\genparam^S}^2 > 1 \text{ or }
\vecnorm{\genparam^S - \genparamstar}{2}^2 >
\frac{4}{1-\enorm{\genparamstar}^2} \delta^2 \Big) &\le \Exs
\enorm{\genparam^S}^2 + \frac{1 - \enorm{\genparamstar}^2}{4 \delta^2}
\Exs \vecnorm{\genparam^S - \genparamstar}{2}^2 \\ &<
\frac{\enorm{\genparamstar}^2 + 1}{2} + \frac{1 -
  \enorm{\genparamstar}^2}{2} = 1.
\end{align}
We conclude that there exists some sign vector $z^S$ satisfying both
inequalities (i) and (iii).


\subsection{Proof of \autoref{LemSudakov}}
\label{SecPfLemSudakov}

We prove \autoref{LemSudakov} via two successive applications of the
Sudakov-Fernique comparison inequality.  In order to keep our
presentation self-contained, let us restate a version of this result
here (e.g., see Theorem 3.15 in \citet{LedTal91}).
For a given a pair of centered Gaussian vectors $\{X_j,~j =
1,\ldots,N\}$ and $\{Y_j,~j = 1,\ldots,N\}$, suppose that
\begin{align*}
  \var(X_i-X_j) \leq \var(Y_i - Y_j) \qquad \mbox{for all $(i,j) \in
    [N]\times [N]$.}
\end{align*}
The Sudakov-Fernique comparison then asserts that $\Exs [\max
  \limits_{j \in [N]} X_j] \leq \Exs[\max \limits_{j \in [N]} Y_j]$.

Using this result, we now prove our claim.  For each $S \in \Tclass$,
define the zero-mean Gaussian random variable $g^S \defn
\inprod{\sgauss}{U \Mat V z^S}$.  First, define a diagonal matrix $D
\defn \diag( 0,\ldots,0, \underbrace{1, \ldots, 1}_{m_1 : m_2}, 0,
\ldots, 0)$, and the zero-mean Gaussian random variables $\gtil^S
\defn \frac{3}{4} \inprod{w}{D z^S}$.  We claim that the
Sudakov-Fernique comparison implies that
\begin{align}
\label{EqnBigStep2}
\Exs \max_{S \in \Tclass} \inprod{\sgauss}{U \Mat V z^S} \ge
\frac{1}{4} \Exs \max_{S \in \Tclass} \inprod{\sgauss}{D z^S}.
\end{align}
See below for the details of this claim.  Second, we introduce the
vector $\zhat^S$ with components $\zhat_i^S \defn |z_i^S|$, and define
a third Gaussian process using the variables $\gring^S \defn
\inprod{\sgauss}{D (\zhat^S - \zhat^{S'})}$.  We also claim that
\begin{align}
\label{EqnBigStep3}
\Exs \max_{S \in \Tclass} \inprod{\sgauss}{D z^S} \ge \Exs
\Big[\max_{S \in \Tclass} \sum_{i \in S} \sgauss_i \Big].
\end{align}
These two claims in conjunction imply the claim of
\autoref{LemSudakov}.  Let us now prove
inequalities~\eqref{EqnBigStep2} and \eqref{EqnBigStep3}.

\paragraph*{Proof of inequality~\eqref{EqnBigStep2}:}
We claim that the processes
$\{g^S, S \in \Tclass\}$ and $\{\gtil^S, S \in \Tclass \}$ satisfy the
Sudakov-Fernique conditions.
In order to prove this claim, we need to verify that for
all subsets $S, S' \in \Tclass$, we have relation
\mbox{$\var(g^S - g^{S'})
  \geq \var(\gtil^S - \gtil^{S'})$.}  On one hand, we have
\begin{align}
\label{EqnSFCondition}
\var(g^S - g^{S'}) \; = \; \Exs\inprod{\sgauss}{U \Mat V (z^S -
  z^{S'})}^2 = \vecnorm{U \Mat V (z^S - z^{S'})}{2}^2 = \vecnorm{H V
  (z^S - z^{S'})}{2}^2,
\end{align}
where the last step uses the orthonormality of $U$.  On the other
hand, we have the equality \mbox{$\var(\gtil^S - \gtil^{S'}) =
  \vecnorm{ D (z^S - z^{S'}) }{2}^2$.}  Consequently, it suffices to
show that there exists an orthogonal matrix $V$ such that
\begin{align}
\label{EqnEigLowerBound}
(\Mat V)^\top \Mat V \succeq \frac{1}{16} D^2.
\end{align}
In order to see this fact, part (e) of \autoref{LemPacking} implies
that the $m_2$ largest eigenvalues of $\Mat^\top \Mat = V \Sigma^2
V^\top$ is lower bounded by $1 - \frac{m_2}{\lwk - 1} -
\sqrt{\frac{a^2 - 9 b^2}{9 b^2}}$.  With the choice of the constants
$(a,b)$ specified above (see paragraph below \autoref{LemPacking}), it
is guaranteed that $\frac{a^2 - 9 b^2}{9b^2} \le \frac{1}{4}$.  This
observation and the definition $m_2 \defn \floor{(\lwk - 1) / 4}$
together imply that
\begin{align}
\label{EqnSVbound}
1 - \frac{m_2}{\lwk - 1} - \sqrt{\frac{a^2 - 9 b^2}{9 b^2}} \ge 1 -
\frac{1}{4} - \frac{1}{2} = \frac{1}{4},
\end{align}
which implies the claim~\eqref{EqnEigLowerBound},
and further completes the proof of the lower bound~\eqref{EqnBigStep2}.

\paragraph*{Proof of inequality~\eqref{EqnBigStep3}:}
The vector $\zhat^S$ defined above is an indicator vector for the
support of $z^S$.  Defining a third Gaussian process using the
variables $\gring^S \defn \inprod{\sgauss}{D (\zhat^S - \zhat^{S'})}$,
we have
  \begin{align}
    \var(\gtil^S - \gtil^{S'}) = \vecnorm{D (z^S - z^{S'})}{2}^2 & \;
    \geq \; \vecnorm{D (\zhat^S - \zhat^{S'})}{2}^2 \; = \;
    \var(\gring^S - \gring^{S'}).
\end{align}
A second application of the Sudakov-Fernique inequality then yields
\begin{align}
\Exs \max_{S \in \Tclass} \inprod{\sgauss}{D z^S} \ge \Exs \max_{S \in
  \Tclass} \inprod{\sgauss}{D \zhat^S} = \Exs \Big[\max_{S \in
    \Tclass} \sum_{i \in S} \sgauss_i \Big],
\end{align}
where in the last step we recall the fact that $S$ is supported on the
set $\{m_1, \ldots, m_2\}$.


\section{Proof of auxiliary lemmas}

In this appendix, we collect the proofs of various auxiliary lemmas.


\subsection{Proof of \autoref{LemKolLinf}}
\label{SecProofLemKolinf}

Let us first state the lemma used in \autoref{SecExtremal}.
\begin{lems}
\label{LemKolLinf}
For an extremal vector of the form $\genparamstar = \sqrt{\eig{s}}e_s
- re_s$, the critical dimension~\eqref{Eqnkcrit} is lower bounded as
\begin{align}
\kind(\delta) \geq 0.09 \cdot \arg \max_{1 \leq k \leq \usedim} \Big
\{\mu^2_k \geq \delta^2 \mu_s \Big \}.
\end{align}
\end{lems}
\noindent The rest of this section is devoted to the proof of this
lemma.

Defining the integer $m \defn \max \{2, \arg \max_{1 \leq k \leq
  \usedim} \Big \{\mu^2_k \geq \delta^2 \mu_s \Big \} \}$,
\citet{wei2017testing} show that we can inscribe an
\mbox{$(m-1)$-dimensional} $\ell_{\infty}$ ball with radius
$\delta/\sqrt{m-1}$ into the ellipse $\ellip$; in particular, see
Section 4.4 of the paper~\cite{wei2017testing}.  We claim that the
Kolmogorov $k$-widths of the $s$-dimensional $\ell_{\infty}$ ball of
radius $\frac{1}{\sqrt{s}}$ are lower bounded as
\begin{align}
\label{EqnKolLinf}
\kwidth_{k}(\Ball_{\infty}(1/\sqrt{s})) \geq 1 - \frac{k}{s}.
\end{align}
Taking this claim as given for the moment, we use it to complete the
proof of \autoref{LemKolLinf}.  Using the lower
bound~\eqref{EqnKolLinf}, we have
\begin{align*}
\kwidth_{k}\Big(\Ellipse_{\at} \cap \Ball(\xi)\Big) \geq
\kwidth_{k}\Big(\Ball_{\infty}(\frac{\xi}{\sqrt{m-1}})\Big) \geq
\left(1 - \frac{k}{m-1} \right)\xi.
\end{align*}
With $\xi \defn (1-\tinyconst) \delta$ and $k = (1 -
\frac{0.9}{1-\tinyconst}) (m - 1)$ the above becomes
$\kwidth_{k}\Big(\Ellipse_{\at} \cap \Ball((1-\eta) \delta)\Big) \ge
0.9 \delta$, so by the definition of the critical
dimension~\eqref{Eqnkcrit}, we have
\begin{align}
\kind(\delta) \geq (1 - \frac{0.9}{1-\eta})(m-1) \geq 0.09 \cdot \arg
\max_{1 \leq k \leq \usedim} \Big \{\mu^2_k \geq \delta^2 \mu_s \Big
\},
\end{align}
as claimed.\\

\noindent The only remaining detail is to prove inequality~\eqref{EqnKolLinf}.


\paragraph*{Proof of inequality~\eqref{EqnKolLinf}:}

Define the set $\Vset \defn \{v \in \real^s \mid v_i =
\pm\frac{1}{\sqrt{s}} \}$ with cardinality \mbox{$M = 2^s$.} We claim
that for any $k$-dimensional subspace $W \subseteq \real^s$, there
exists some $v \in \Vset$ such that
\begin{align}
  \label{EqnRabbit}
\ltwo{v - \Pi_W(v)}^2 \geq 1 - \frac{k}{s}.
\end{align}
Then by definition of Kolmogorov width, the
inequality~\eqref{EqnKolLinf} holds.  In order to prove the lower
bound~\eqref{EqnRabbit}, we take an orthonormal basis $z_1,\ldots,z_k$
of $W$ and extend it to an orthonormal basis $z_1,\ldots,z_m$ for
$\real^s$. We then have
\begin{align*}
\sum_{v \in \Vset} \ltwo{v - \Pi_W(v)}^2 = \sum_{v \in \Vset} \sum_{j=k+1}^s
\inprod{v}{z_j}^2 = \sum_{j=k+1}^s \sum_{v \in \Vset} \inprod{v}{z_j}^2
= (s-k) \cdot \frac{M}{s},
\end{align*}
where we have used the fact that
\begin{align*}
    \sum_{v \in \Vset} \inprod{v}{z_j}^2 = M \cdot \frac{1}{s}
    \ltwo{z_j}^2 = \frac{M}{s}.
\end{align*}
Therefore, there must exist some $v \in \Vset$ such that $\ltwo{v -
  \Pi_W(v)}^2 \geq 1-\frac{k}{s}$, which establishes the
inequality~\eqref{EqnKolLinf}.


\subsection{Proof of \autoref{Lemfun}}
\label{SecProofLemfun}

Recalling our calculations from \autoref{ExampleSobolev}, we found
that
\begin{align*}
\kind(0, \frac{1}{2} \delta) = \argmin_k \{\sqrt{\eig{k+1}} \le
\frac{9}{10} \cdot \frac{1}{2} \delta\}.
 \end{align*}
Consequently, in order to prove \autoref{Lemfun}, it suffices to
show that
\begin{align}
\label{EqnfunKey}
\upk(\genparamstar, \delta) \leq \arg \min_k \{2\sqrt{\mu_k} \leq
\frac{9}{10}\delta\} \qquad \mbox{for all }\delta \leq \sqrt{\mu_1}.
\end{align}
By definition of the critical dimension~\eqref{Eqnkcrit}, it is
sufficient to show that the Kolmogorov width is upper bounded as
\begin{align}
\label{EqnEasy}
\kwidth_k(\Ellipse_{\at} \cap \Ball(a \delta)) \leq
\min\{a \delta, ~2\sqrt{\mu_k}\},
\end{align}
where $a \defn 1 - \tinyconst$.

We claim that the set $\Ellipse_{\at} \cap \Ball(a \delta)$ is
contained within the set $2\Ellipse \cap \Ball(a \delta)$.  Indeed,
note that any $v \in \Ellipse_{\at} \cap \Ball(a \delta)$ has
Euclidean norm bounded as $\ltwo{v}\leq a\delta$ and Hilbert norm
bounded as $\enorm{v + \at} \leq 1$. The Cauchy-Schwarz further
guarantees that
\begin{align*}
\enorm{v}^2 = \enorm{v + \at - \at}^2 \leq 2 \enorm{v + \at}^2 + 2
\enorm{\at}^2 \leq 4,
\end{align*}
where the last step follows from the fact that both $\at$ and $v +
\at$ lie in ellipse $\Ellipse$.  We have thus established the claimed
set inclusion.

From this set inclusion, we have
\begin{align*}
\kwidth_k(\Ellipse_{\at} \cap \Ball(a \delta)) \leq
\kwidth_k(2\Ellipse \cap \Ball(a \delta)) = \min\{a \delta,
~2\sqrt{\mu_k}\},
\end{align*}
which establishes the claim~\eqref{EqnEasy}.  Putting pieces together
completes the proof of \autoref{Lemfun}.


\section{Well-definedness of the function \texorpdfstring{$\Rfun$}{Phi}}
\label{AppRfun}

\begin{figure}[h]
  \centering \includegraphics[width=0.45\textwidth]{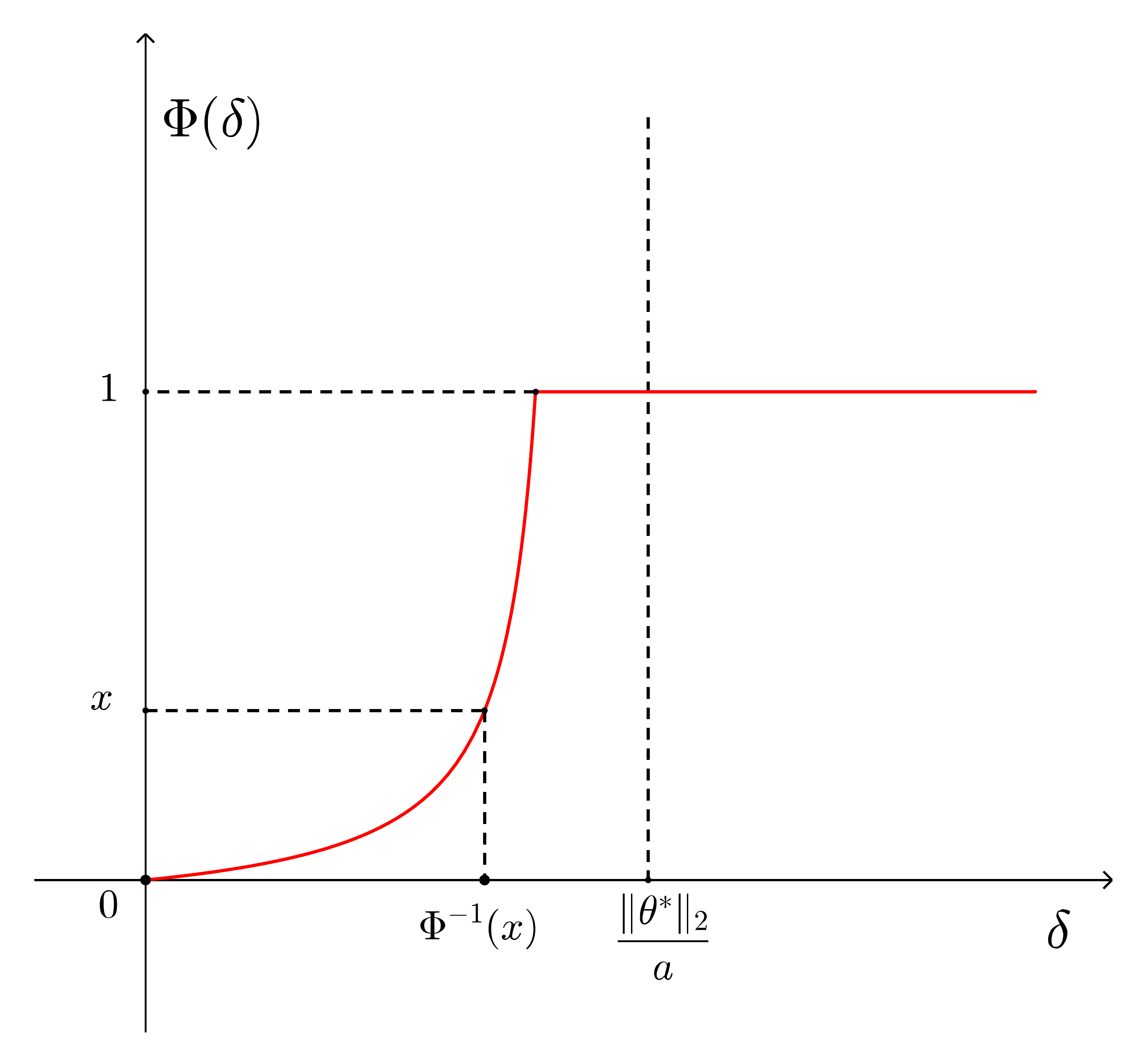}
  \caption{Illustration of the function $\Rfun$.}
  \label{fig:Rfun}
\end{figure}

In this appendix, we verify that the function $\Rfun$ from
equation~\eqref{EqnRfun} is well-defined.
We again use the shorthand $a \defn 1 - \tinyconst$.
In order to provide
intuition, \autoref{fig:Rfun} provides an illustration of $\Rfun$.

We begin with the case when $\delta < \ltwo{\at}/a$.  For simplicity
of notation, let
\begin{align}
r(\delta) \defn \min \Big \{ r \geq 0 \, \mid a^2 \delta^2 \leq
\sum_{i=1}^\usedim \frac{r^2}{(r+\mu_i)^2} (\at_i)^2 \Big \}.
\end{align}
Note that for each for each $\mu_i \geq 0$, the function $f(r) \defn
\sum_{i=1}^\usedim \frac{r^2}{(r+\mu_i)^2}(\at_i)^2$ is non-decreasing
in $r$.  It is also easy to check that
\begin{align*}
 \lim_{r\to 0^+} f(r) = 0,\qquad \text{ and }\lim_{r\to \infty} f(r) =
 \ltwo{\at}^2.
\end{align*}
Then the quantity $r(\delta)$ is uniquely defined and positive
whenever $\delta < \ltwo{\at}/a$.  Note that as $\delta \to
\frac{\ltwo{\at}}{a}$, $a^2 \delta^2 \to \ltwo{\at}^2$ therefore
$r(\delta) \to \infty$.

It is worth noticing that given any $\at$ where $\ltwo{\at}$ does not
depend on $\delta$, $r$ goes to zero when $\delta \to 0$, namely
$\lim_{\delta \to 0^+} \Rfun(\delta) = 0$.

\vspace*{2cm}




\end{document}